\documentclass[a4paper,12pt]{amsart}
\usepackage{a4wide,amssymb}
\usepackage{xypic}

\title [Generalised Hopf type formulas]{\bf N-fold \v{C}ech Derived Functors
and Generalised Hopf type formulas}
\author[Guram Donadze, Nick Inassaridze and Timothy Porter]{}
\keywords{\v{C}ech derived functors, Hopf type formulas, $K$-functors}
\subjclass{18G50, 18G10}

\newtheorem{theorem}{Theorem}
\newtheorem{proposition}[theorem]{Proposition}
\newtheorem{lemma}[theorem]{Lemma}
\newtheorem{corollary}[theorem]{Corollary}

\def\two{\begin{matrix}\to\\[-3mm]\to\end{matrix}}
\def\atwo#1#2{\begin{matrix}&\ovs{#1}{\to}\\[-3mm]&\uns{#2}{\to} \end{matrix}}
\def\three{\begin{matrix}\to\\[-2mm]\to\\[-2mm]\to\end{matrix}}
\def\athree#1#2{\begin{matrix}\ovs{#1}{\to}\\[-2mm]\to\\[-2mm]\uns{#2}{\to}\end{matrix}}

\def\dtwo{\begin{matrix}\to\\[-3.5mm]\vdots\\[-3mm]\to\end{matrix}}
\def\adtwo#1#2{\begin{matrix}&\ovs{#1}{\to}\\[-3.5mm]&\vdots\\[-3mm]&\uns{#2}{\to}\end{matrix}}

\def\al{\alpha}
\def\be{\beta}
\def\ka{\kappa}
\def\de{\delta}
\def\De{\Delta}

\def\vp{\varphi}
\def\ep{\epsilon}

\def\pa{\partial}

\def\lam{\lambda}
\def\Gm{\Gamma}

\def\lra{\longrightarrow}
\def\lla{\longleftarrow}
\def\hra{\hookrightarrow}
\def\la{\langle}
\def\ra{\rangle}
\def\da{\downarrow}

\def\lbr{\linebreak}

\def\ol{\overline}
\def\ul{\underline}
\def\uns{\underset}
\def\ovs{\overset}
\def\stmin{\setminus}
\def\seq{\subseteq}
\def\esq{\supseteq}
\def\eset{\emptyset}
\def\wt{\widetilde}

\def\Ker{\operatorname{Ker}}
\def\Im{\operatorname{Im}}

\def\ob{\operatorname{ob}}
\def\cat{\operatorname{cat}}
\def\lim{\operatorname{lim}}
\def\min{\operatorname{min}}

\advance\headheight1.15pt

\begin{document}

\maketitle

\noindent
{\bf Guram Donadze\\
{\footnotesize \it A.Razmadze Mathematical Institute, Georgian
Academy of Sciences \\ M.Alexidze St. 1, Tbilisi 380093. Georgia
\\ donad@rmi.acnet.ge}}

\noindent
{\bf Nick Inassaridze\\
{\footnotesize \it A.Razmadze Mathematical Institute, Georgian Academy of
Sciences \\ M.Alexidze St. 1, Tbilisi 380093. Georgia \\
inas@rmi.acnet.ge}}

\noindent
{\bf Timothy Porter\\
{\footnotesize \it Mathematics Division, School of Informatics, University of Wales
Bangor,\\
Bangor, Gwynedd LL57 1UT,  United Kingdom.\\
t.porter@bangor.ac.uk}}

\

\begin{abstract}
In 1988, Brown and Ellis published \cite{be} a generalised Hopf
formula for the higher homology of a group. Although substantially
correct, their result lacks one necessary condition. We give here
a counterexample to the result without that condition. The main
aim of this paper is, however, to generalise this corrected result
to derive formulae of Hopf type for the $n$-fold \v{C}ech derived
functors of the lower central series functors $Z_k$. The paper
ends with an application to algebraic K-theory.
\end{abstract}

\

\section*{Introduction and Summary}

\

The well known Hopf formula for the second integral homology of a
group says that for a given group $G$ there is an isomorphism
$$
H_2(G)\cong \frac{R\cap [F,F]}{[F,R]}\;,
$$
where $R\rightarrowtail F\twoheadrightarrow G$ is a free presentation of the group $G$.

Several alternative generalisations
of this classical Hopf formula to higher dimensions were made in various papers,
\cite{con, Ro, St}, but perhaps the most successful one, giving formulas in
all dimensions, was by Brown and Ellis, \cite{be}.
They used topological methods,  and in particular the Hurewicz theorem
for $n$-cubes of spaces, \cite{bl2}, which itself is an
application of the generalised van Kampen theorem for diagrams of
spaces \cite{bl}.  The end result was:

{\noindent \bf Theorem BE (\cite{be}).} {\em Let $R_1$, \ldots,
$R_n$ be normal subgroups of a group $F$ such that
$$
\begin{matrix}
H_2(F)=0,\;\;\; H_r(F\diagup \underset{i\in
A}{\prod}R_i)=0,\;\;\;\text{for}\;\; r=|A|+1,\; r=|A|+2,\;
\end{matrix}
$$
with $A$ a non-empty proper subset of $\la n\ra =\{ 1,\ldots ,n\}$
(for example, if the groups $F\diagup \uns{i\in A}{\prod}R_i$ are
free for $A\neq\la n\ra$) and $F\diagup \underset{1\le i\le n}
{\prod}R_i\cong G$. Then there is an isomorphism
$$
H_{n+1}(G)\cong \frac{ \overset{n}{\underset{i=1}{\cap}}R_i\cap
[F,F]} {\underset{A\seq \la n \ra}{\prod} [\underset{i\in
A}{\cap}R_i, \underset{i\notin A}{\cap}R_i]}\;.
$$}Later, Ellis using mainly algebraic means, and,  in particular,
his hyper-relative derived functors, proved the same result,
\cite{e}.

The similarity between this formula and the formulae given by Mutlu and the
third author for the various homotopy invariants of a simplicial group (see
\cite{atmp1,atmp2}) suggested that there should be a purely algebraic proof of
this, which hopefully would generalise further.  Examining the classical case ($n =
1$), and the proof of the usual Hopf formula, showed a link with the \v{C}ech
derived functors of the abelianisation functor, (cf. \cite{HI}).

Trying to derive this result purely algebraically and to obtain
Hopf type formulas for some more general situations, we suspected that
the conditions given above for Theorem BE were not sufficient for
getting the generalised Hopf formula for $H_{n+1}(G)$, $n\ge 3$.
In fact, we give the following counter-example to Theorem BE:

Let $F$ be a free group with base $\{x_1,x_2\}$, $R_1$, $R_2$ and
$R_3$ normal subgroups of the group $F$ generated by the one point
sets $\{x_1\}$,  $\{x_2\}$ and $\{x_1x_2^{-1}\}$ respectively and
$G=1$. Then we have $F/R_i\cong {\mathbb Z}$, $i=1,2,3$, $F/R_i
R_j=1$, $i\neq j$ and $[F,F]=[R_1,R_2]=R_i\cap R_j$, $i\neq j$,
therefore $$\frac{\uns{i\in \la 3\ra}{\cap}R_i\cap
[F,F]}{\uns{A\seq \la 3\ra}{\prod}[\uns{i\in
A}{\cap}R_i,\uns{i\notin A}{\cap}R_i]}\cong {\mathbb Z}$$ whilst
$H_3(G)=1$.

We thus set out to prove a corrected version of this Brown-Ellis
generalised Hopf formulae, but also to generalise it further in
the following direction. Homology groups are the derived functors
of the abelianization functor. Our generalisation handles the
derived functors of the functors that kill higher commutators.
More precisely, let the endofunctors $Z_k(G)$ be given by
$Z_k(G)=G/\Gm_k(G)$, $k\ge 2$, where $\{\Gm_k(G),\; k\ge 1\}$ is
the lower central series of a given group $G$.  These $Z_k$ are
endofunctors on the category of groups and  generalise the
abelianization functor, so their non-abelian left derived
functors, $L_nZ_k$, $n\ge 0$, generalise the group homology
functors $H_n$, $n\ge 1$, cf., for instance, \cite{bb}.

In \cite{IE}, a  Hopf-like formula is proved for the second
Conduch\'{e}-Ellis homology of precrossed modules using
\emph{\v{C}ech derived functors}. The main goal of this paper is
to develop this method further, and by applying it, to  express
$L_nZ_k$, $n\ge 1$, $k\ge 2$, by generalised Hopf type formulas.
In particular, we will give the stronger conditions needed for
Theorem BE. Finally we apply these results to algebraic
$K$-theory.

\medskip

In the first section, we introduce the \v{C}ech derived functors
illustrating their use by proving the classical Hopf formula in a
new way. This is not just an illustration as it does indicate some
of the ways the argument will go later on.

In Section 2 we introduce the notion of \emph{simple} normal
$(n+1)$-ad of groups $(F;R_1,\ldots, R_n)$ \emph{relative to}
$R_j$ for some $1\le j \le n$. Then we show that for a given
pseudosimplicial group $F_*$, the normal $(j+1)$-ads of groups
$(F_n;\Ker d^n_0,\ldots,\Ker d^n_{j-1})$ are simple relative to
$\Ker d^n_{j-1}$ for all $1\le j\le n$ (Proposition \ref{3}). The
main result of Section 2 is Theorem \ref{5}, giving that for an
aspherical augmented pseudosimplicial group $(F_*,d^0_0,G)$, there
is a natural isomorphism
$$\pi_nZ_k(F_*)\cong \frac {\ovs{n-1}{\uns{i=0}{\cap}}\Ker
{d^{n-1}_i}\cap \Gm_k(F_{n-1})}{D_k (F_{n-1};\Ker
d^{n-1}_0,\ldots,\Ker d^{n-1}_{n-1})}$$ for $n\ge 1$, $k\ge 2$.
Here the $D_k$-term takes the form of an iterated commutator
subgroup which is an obvious generalisation of the denominator
terms of both the classical Hopf formula and the Brown-Ellis
extension of that formula.  It is also related to the descriptions
of the image term of the Moore complex, as used  in
\cite{atmp1,atmp2}.  The explicit formula is given at the start of
Section 2.

For an inclusion crossed $n$-cube of groups, ${\mathcal M}$, given
by a normal $(n+1)$-ad of groups we construct, in Section 3, a new
induced crossed $n$-cube ${\mathcal B}_k({\mathcal M})$, $k\ge 2$
(Proposition \ref{7}). We show the existence of an isomorphism of
simplicial groups $Z_kE^{(n)}({\mathcal M})_*\cong
E^{(n)}({\mathcal B}_k({\mathcal M}))_*$, where $E^{(n)}({\mathcal
M})_*$ denotes the diagonal of the $n$-simplicial nerve of the
crossed $n$-cube of groups ${\mathcal M}$, (Proposition \ref{8}).

Section 4 is devoted to the investigation of some properties of
the mapping cone complex of a morphism of (non-abelian) group
complexes introduced in \cite{l1}. In particular, for a given
morphism of pseudosimplicial groups $\al:G_*\to H_*$ the natural
morphism $\ka: NM_*(\al)\to C_*(\wt\al)$ induces isomorphisms of
their homology groups, where $C_*(\wt\al)$ is the mapping cone
complex of the induced morphism of the Moore complexes and
$NM_*(\al)$ is the Moore complex of a new pseudosimplicial group
constructed using $\al$ (Proposition \ref{9}). (Here similar
results have recently been found by Conduch\'e, \cite{conduche2}.)
Using this result we derive  purely algebraicly the result of
\cite{l1}, (3.4. Proposition), giving for a crossed $n$-cube of
groups ${\mathcal M}$ an isomorphism between the homotopy groups
of $E^{(n)}({\mathcal M})_*$ and the corresponding homology groups
of a chain complex of groups $C_*({\mathcal M})$, (Proposition
\ref{10}). In particular, we give an explicit computation of  the
$n$-th homotopy group of the simplicial group $E^{(n)}({\mathcal
M})_*$.

In Section 5, we introduce a notion of $n$-fold \v{C}ech derived
functors of an endofunctor on the category of groups (Theorem
\ref{12}, Definition) which will be the subject of future papers
and applications to nonabelian homological algebra and $K$-theory.
We give an explicit calculation of $n$-th $n$-fold \v{C}ech
derived functor of the functor $Z_k$, $k\ge 2$ (Theorem \ref{16}).
Our method gives the  possibility of finding the sought for
sufficient conditions for,  and a purely algebraic proof of,  the
generalised Hopf formula of Brown and Ellis, moreover we express
$L_nZ_k(G)$, $n\ge 1$, $k\ge 2$ by a Hopf type formula (Theorem
\ref{17}).

In Section 6, an application to algebraic $K$-theory is given. In
particular, Quillen's algebraic  $K$-functors $K_{n+1}$, $n\ge 1$
are described in terms of short exact sequence including the
higher Hopf type formulae for free exact $n$-presentations induced
by a free simplicial resolution of the general linear group
(Theorem \ref{21}).

\

%%Section 1
\section{An approach to the classical Hopf formula via  \v{C}ech derived functors}

\

We give here a brief introduction to \v{C}ech derived functors.  A
fuller account is given in \cite{HI}.  We will limit ourselves to
the \v{C}ech derived functors of the Abelianization functor.
Later we will develop the $n$-fold analogue of some of this
theory.

\medskip

{\noindent \em Definition.} Let $T:{\mathfrak Gr}\to {\mathfrak
Gr}$ be a covariant functor. Define \emph{$i$-th  \v{C}ech derived
functor} ${\mathcal L}_iT:{\mathfrak Gr}\to {\mathfrak Gr}$, $i\ge
0$, of the functor $T$ by choosing for each $G$ in ${\mathfrak
Gr}$, a free presentation ${\mathfrak F} :  R\rightarrowtail
F\stackrel\al\twoheadrightarrow G$ of $G$ and setting
$$
{\mathcal L}_iT(G)= \pi_i(T\check{C}(\al)_*)\;,
$$
where $(\check{C}(\al)_*,\al,G)$ is the \v{C}ech resolution of the
group $G$ for the free presentation ${\mathfrak F}$ of $G$. This
latter resolution is constructed as follows:

Given a group $G$ and a homomorphism of groups $\al:F\to G$. The
\emph{\v{C}ech augmented complex} $(\check{C}(\al)_*,\al,G)$ for
$\al$ is
$$
\begin{matrix}
 \cdots \dtwo & F\times_G \cdots \times_G F & \adtwo{d^n_0}{d^n_n} \cdots
 \athree{d^2_0}{d^2_2}
& F\times_G F & \atwo{d^1_0}{d^1_1} & F & \ovs{\al}{\lra} & G\;,
\end{matrix}
$$
thus
$$
\check{C}(\al)_n=\underbrace{F\times_G \cdots \times_G
  F}_{(n+1)-\text{times}}=\{(x_0,\ldots,x_n)\in F^{n+1}\;|\;\al(x_0)=\cdots =\al(x_n)\}\;\;\text{for}\;\;n\ge 0\;,
$$
$$
d^n_i(x_0,\ldots,x_n)=(x_0, \ldots, \hat{x_i}, \ldots, x_n)
$$
and
$$
s^n_i(x_0,\ldots,x_n)=(x_0, \ldots, x_i, x_i, x_{i+1},\ldots, x_n)
$$
for $0\le i \le n$ (see \cite{HI}).

In case $F$ is a free group and $\al$ is an epimorphism as in
$\mathfrak{F}$ above, $(\check{C}(\al)_*,\al,G)$ will be called a
\emph{\v{C}ech resolution} of $G$.
\medskip
\emph{Example. }   The prime example of the \v{C}ech derived functors are those of the abelianization functor.  We recall that
$$
H_2(G)\cong {\mathcal L}_1{\it Ab}(G)\;$$and will use this later.
\medskip

Now using crossed modules and their nerves, we present a fresh
view of the \v{C}ech complex which leads to some  ideas that will
be useful throughout the paper.

First, let us recall the definition of crossed module. A crossed
module $(M,P,\mu)$ is a group homomorphism $\mu:M\to P$ together
with a (left) action of $P$ on $M$ which satisfies the following
conditions:
\begin{enumerate}
\item[(i)] $\mu(^pm)=p\mu(m)p^{-1}$,
\item[(ii)]  $^{\mu(m)}m'=mm'm^{-1}$, \qquad \quad(Peiffer identity)
\end{enumerate}
for all $m,m'\in M$ and $p\in P$.

A morphism $(\vp,\psi):(M,P,\mu)\to (N,Q,\nu)$ of crossed modules
is a commutative square of groups
$$
\xymatrix{
M \ar[r]^{\vp}\ar[d]_\mu& N \ar[d]^\nu\\
P \ar[r]^{\psi}& Q}
$$
such that $\vp(^pm)= {}^{\psi(p)}{\vp(m)}$ for all $m\in M$, $p\in
P$. Let us denote the category of crossed modules by ${\mathcal{
CM}}$.

It is well known from \cite{l1} that the category ${\mathcal CM}$
of crossed modules is equivalent to the category of ${\it
cat}^1$-groups (for the definition see \cite{l1}), and for a given
crossed module, ${\mathcal M}=(M\stackrel{\mu}{\rightarrow}P)$,
the corresponding ${\it cat}^1$-group is $(M\rtimes P, s, t)$,
where $s(m,p)=p$ and $t(m,p)=\mu(m)p$. This ${\it cat}^1$-group
has an internal category stucture within the category
$\mathfrak{Gr}$ of groups and the nerve of its category structure
forms the following simplicial group
$$
E({\mathcal M})_*
:\quad\xymatrix{\ldots\ar[r]<1.6ex>\ar[r]<-1.6ex>^{\vdots\hspace{.5cm}}
& E({\mathcal
M})_n\ar[r]<1.6ex>\ar[r]<-1.6ex>^{\hspace{.5cm}\vdots} &
\quad\cdots\quad\ar[r]<1ex>^{d_0}\ar[r]<-1ex>_{d_2} \ar[r]&
E({\mathcal M})_1\ar[r]<.5ex>^{d_0}\ar[r]<-.5ex>_{d_1} &
E({\mathcal M})_0 ,}
$$
where $E({\mathcal M})_n = M\rtimes (\cdots (M\rtimes P)\cdots )$
with $n$ semidirect factors of $M$, and face and degeneracy
homomorphisms are defined by
\begin{align*}
&d_0(m_1,\ldots,m_n,p)=(m_2,\ldots,m_n,p)\;,\\
&d_i(m_1,\ldots,m_n,p)=(m_1,\ldots,m_im_{i+1},\ldots,m_n,p)\;,\;\;\;0<i<n\;,\\
&d_n(m_1,\ldots,m_n,p)=(m_1,\ldots,m_{n-1},\mu(m_n)p)\;,\\
&s_i(m_1,\ldots,m_n,p)=(m_1,\ldots,m_i,1,m_{i+1},\ldots,m_n,p)\;,\;\;\;0\le
i\le n\;.
\end{align*}
\noindent The simplicial group $E({\mathcal M})_*$ is called the
\emph{nerve} of the crossed module ${\mathcal M}$ and its Moore
complex is trivial in dimensions $\ge$ 2. In fact its Moore
complex is just the original crossed module up to isomorphism with
$M$ in dimension 1 and $P$ in dimension 0.

\medskip

%%Lemma 1
\begin{lemma}\label{1.1}
Let $G$ be a group and $F\stackrel{\al}{\rightarrow}G$ be a
homomorphism of groups. Then the \v{C}ech complex for $\al$ is
isomorphic to the nerve of the inclusion crossed module $E(R\hra
F)_*$, where $R = \Ker \al$.
\end{lemma}
\begin{proof} Let us compare the following two simplicial groups, the
\v{C}ech complex for $\al$ and the nerve of the inclusion crossed
module $R\hra F$
$$
\xymatrix{E(R\hra F)_* \hspace{-1cm}&\quad
:\quad\ldots\ar[r]<1.6ex>\ar[r]<-1.6ex>^{\vdots\hspace{.5cm}} &
  R\rtimes R\rtimes F\ar[r]<1ex>\ar[r]<-1ex> \ar[r]\ar[d]^{\lam_2}&R\rtimes
  F\ar[r]<.5ex>\ar[r]<-.5ex>\ar[d]^{\lam_1} & F\ar[d]^{\lam_0} \\
\check{C}(\al)_* \hspace{-1cm}&\quad
:\quad\ldots\ar[r]<1.6ex>\ar[r]<-1.6ex>^{\vdots\hspace{.5cm}} &
  F\times_G F\times_G F\ar[r]<1ex>\ar[r]<-1ex> \ar[r]&F\times_G F\ar[r]<.5ex>\ar[r]<-.5ex> & F}
$$
\medskip
\noindent by constructing a morphism $\lam_*$ as follows: \\
\\
\hspace*{1cm}$\lam_0$ is the identity on $F$; and  \\
\hspace*{1cm}$\lam_n(r_1,\ldots,r_n,f)=(r_1\cdots r_n f, r_2\cdots
r_n f, \ldots,r_n f, f)$ {} for all {} $n\ge 1$ {} and \\
\hspace*{4.3cm} $(r_1,\ldots,r_n,f)\in E(R\hra F)_n$.

It is easy to check that $\lam_*$ is isomorphism of simplicial
groups.
\end{proof}

\medskip

We recall from \cite{LG} that a crossed module $\mu : M \to P$ is
called \emph{abelian} if $P$ is an abelian group and the action of
$P$ on $M$ is trivial. This implies that $M$ is also abelian. Let
us denote the category of abelian crossed modules by ${\mathfrak
Ab}{\mathcal CM}$.

One can define the abelianization functor $\it Ab$ from the
category ${\mathcal CM}$ to the category ${\mathfrak Ab}{\mathcal
CM}$ in the following way: for any crossed module ${\mathcal M}=(M
\stackrel{\mu}{\to}P)$,
$$
{Ab}({\mathcal M})=(\frac{M}{ [P,M]}\stackrel{\ol\mu}{\to}\frac{P}{[P,P]})\;,
$$
in which $[P,M]$ is the subgroup of $M$ generated by the elements
${}^pmm^{-1}$ for all $m\in M$, $p\in P$ and $\ol\mu$ is induced
by $\mu$.

Given a simplicial group $G_*$, let us apply the group
abelianization functor dimension-wise, denote the resulting
simplicial group by ${\it Ab}(G_*)$.

\medskip

%%Proposition 2
\begin{proposition}\label{nervexmod}
Let $M \stackrel{\mu}{\to} P$ be a crossed module. Then there is
an isomorphism of simplicial groups
$$
{\it Ab}(E(M \stackrel{\mu}{\to} P)_*)\cong E({\it
Ab}(M\stackrel{\mu}{\to} P))_*\;.
$$
\end{proposition}
\begin{proof} Let us consider the two simplicial groups
$$
 Ab(E(M\stackrel{\mu}{\to} P)_*) :\xymatrix{\quad
:\quad\ldots\ar[r]<1.6ex>\ar[r]<-1.6ex>^{\vdots\hspace{.5cm}} &
  (M\rtimes M\rtimes P)^{ab}\ar[r]<1ex>\ar[r]<-1ex> \ar[r]&(M\rtimes
  P)^{ab}\ar[r]<.5ex>\ar[r]<-.5ex> & P^{ab} }
$$
and
$$
E({\it Ab}(M\stackrel{\mu}{\to} P))_*:\xymatrix{\quad
:\quad\ldots\ar[r]<1.6ex>\ar[r]<-1.6ex>^{\vdots\hspace{.5cm}} &
  \displaystyle\frac{M}{[P,M]}\times\frac{M}{[P,M]} \times P^{ab}\ar[r]<1ex>\ar[r]<-1ex> \ar[r]&\displaystyle\frac{M}{[P,M]}\times
  P^{ab}\ar[r]<.5ex>\ar[r]<-.5ex> & P^{ab} }.
$$
It is easy to show that $$(M\rtimes P)^{\it
ab}\stackrel{\ka}{\to}\frac {M}{[P,M]}\times P^{\it ab},$$
$$\ka[(m,p)]=([m], [p])$$ and $$\frac {M}{[P,M]}\times P^{\it
ab}\stackrel {\ka'}{\to} (M\rtimes P)^{\it ab},$$ $$\ka'([m], [p])
= [(m,p)]$$ are homomorphisms and $\ka\ka' =\ka'\ka =1$. Using
these isomorphisms one has
\begin{eqnarray*}
(M\rtimes (M\rtimes P))^{\it ab}& \cong& \frac{M}{[M\rtimes P, M]}\times
(M\rtimes P)^{\it ab} \\
 &\cong& \frac{M}{[M\rtimes P,M]}\times \frac{M}{[P,M]}\times  P^{\it ab}\\
&\cong& \frac{M}{[P,M]}\times \frac{M}{[P,M]}\times  P^{\it ab}\;,
\end{eqnarray*}
since $[M\rtimes P,M]= [P,M]$ as
\begin{eqnarray*}
{}^{(m,p)}{m'm^{\prime -1}} &=& m {}^pm' m^{-1}
m^{\prime -1}\\&=&\lbr m {}^pm' m^{\prime -1} m^{-1} mm'm^{-1}m^{\prime -1}\\&=& m({}^pm'
m^{\prime -1})m^{-1}({}^{\mu(m)}{m'}m^{\prime -1}) \in [P,M]
\end{eqnarray*} for all $m,m'\in M$, $p\in P$. It is also easy to see  by the same
sort of argument that  there exist isomorphisms between higher
terms of these simplicial groups and that these isomorphisms are
compatible with face and degeneracy maps.
\end{proof}

\medskip

Now we are ready to revisit the classical Hopf formula.

One can prove the Hopf formula in many ways, but for our later
generalization we will prove it using the \v{C}ech derived
functors.

\medskip

%%Theorem 3
\begin{theorem}\label{1.3}
Let $G$ be a group and $R\stackrel{i}{\rightarrowtail}
F\twoheadrightarrow G$ be a free presentation of the group $G$.
Then there is an isomorphism of groups
$$
H_2(G)\cong \frac{ R\cap[F,F]}{[F,R]}\;.
$$
\end{theorem}
\begin{proof} Using \cite{Pi1, Pi3} (see also Theorem 2.39(ii),
\cite{HI}), there is an isomorphism
$$
H_2(G)\cong {\mathcal L}_1{\it Ab}(G)\;,
$$
where ${\mathcal L}_1{\it Ab}$ is the first \v{C}ech derived
functor of the group abelianization functor.

Now Lemma \ref{1.1} and Proposition \ref{nervexmod} implies an
isomorphism
$$
H_2(G)\cong \pi_1 (E({\it Ab}(R\stackrel{i}{\hra} F))_*)\;.
$$

As we already mentioned above, it is clear that $\pi_1 (E({\it
Ab}(R\hra F))_*)$ is isomorphic to the kernel of the crossed
module
$$
{\it Ab}(R\hra F)= \left(\frac{R}{ [F,R]}\stackrel{\ol i}{\lra}\frac{F}{[F,F]}\right)\;
$$
giving the desired result.
\end{proof}

\

%%Corollary 4
\begin{corollary}\label{1.4}
Let $R\rightarrowtail K\stackrel{\al}{\twoheadrightarrow} G$ be a
presentation of a group $G$ and $H_2(K)=0$. Then there is an
isomorphism of groups
$$
H_2(G)\cong \frac{R\cap [K,K]}{[K,R]}\;.
$$
\end{corollary}
\begin{proof} Consider a free presentation $R''\rightarrowtail F\stackrel{\be}{\twoheadrightarrow} K$ of the
group $K$. Hence, one has also a free presentation
$R'\rightarrowtail F\stackrel{\al\be}{\twoheadrightarrow} G$ of
the group $G$. It is easy to check that there is an exact sequence
of groups
$$
\frac{R''\cap [F,F]}{[F,R'']}\lra \frac{R'\cap [F,F]}{[F,R']}\lra
\frac{R}{[K,R]}\lra \frac{K}{[K,K]} \lra \frac{G}{[G,G]}\lra
1\;\;,
$$
Thus by Theorem \ref{1.3} and the condition that $H_2(K)=0$, one
has the following exact sequence of groups
$$
0\lra H_2(G)\lra \frac{R}{[K,R]}\lra K^{\it ab}\lra G^{\it ab}\lra 1\;\;,
$$
which completes the proof.
\end{proof}

\

%%Section 2
\section{Simple normal $(n+1)$-ad of groups}

\

One of the tools we will be using later is the theory of crossed
$n$-cubes of groups.  These generalise normal $(n+1)$-ads of
groups in the same way that crossed modules generalise normal
subgroups.  We therefore start by developing some techniques for
handling $(n+1)$-ads of groups, relating them to iterated
commutators.

Given a group $F$ and  $n$ normal subgroups,  $R_1,\ldots, R_n$,
then $(F;R_1,\ldots,R_n)$ will be called a \emph{normal $(n+1)$-ad
of groups}. A normal $(n+1)$-ad of groups $(F; R_1, \ldots, R_n)$
is called \emph{simple relative to $R_j$} for some $1\le j\le n$
if there exists a subgroup $F'$ of the group $F$ such that
$$
\begin{matrix}
F'\cap R_j =1\;,\qquad \uns{i\in A}{\cap}R_i = (\uns{i\in
A}{\cap}R_i \cap F')(\uns{i\in A}{\cap}R_i \cap R_j)
\end{matrix}
$$
for all $A\seq \la n\ra\stmin \{j\}$.

For a given $(n+1)$-ad of groups $(F; R_1, \ldots, R_n)$, $A\seq
\la n\ra$ and $k\ge 1$ denote by $D_k(F;A)$ the following normal
subgroup of the group $F$
$$
\uns{A_1\cup A_2\cup \cdots \cup A_k=A}{\prod}[\uns{i\in
A_1}{\cap}R_i,[\uns{i\in A_2}{\cap}R_i, ~\ldots~,[\uns{i\in
A_{k-1}}{\cap}R_i,\uns{i\in A_k}{\cap}R_i] \ldots ]]\;.
$$
Sometimes we write $D_k(F;R_1,\ldots,R_n)$ instead of the notation
$D_k(F;\la n\ra)$ .

\medskip

%%Lemma 5
\begin{lemma}\label{1}
Let $(F; R_1, \ldots, R_n)$ be a normal $(n+1)$-ad of groups which
is simple relative to $R_j$, $1\le j\le n$ and $k\ge 1$, then
$$
D_k(F; A)=(D_k(F; A)\cap F')D_k(F; A\cup \{j\})
$$
for all $A\seq \la n\ra\stmin \{j\}$.
\end{lemma}
\begin{proof} We use induction on $k$. Let $k=1$, then
$$
D_1(F; A)=\uns{i\in A}{\cap}R_i=(\uns{i\in A}{\cap}R_i \cap F')
(\uns{i\in A}{\cap}R_i \cap R_j)= (\uns{i\in A}{\cap}R_i \cap
F')D_1(F; A\cup \{j\})
$$
for $A\seq \la n\ra\stmin \{j\}$.

Proceeding by induction, we suppose that the assertion is true for
$k-1$ and we will prove it for $k$.

The inclusion $(D_k(F; A)\cap F')D_k(F; A\cup \{j\})\seq D_k(F;
A)$ is obvious. It is easy to see that a generator of $D_k(F;A)$
has the form $[x,w]$, where $x\in \uns{i\in B}{\cap}R_i$, $w\in
D_{k-1}(F;C)$, $B,C \seq A\seq \la n\ra\stmin \{j\}$ and $B\cup
C=A$. There exist elements $y\in \uns{i\in B}{\cap}R_i\cap F'$ and
$z\in \uns{i\in B}{\cap} R_i\cap R_j$ such that $x=yz$. One has
$$
[x,w]=[yz,w]=y[z,w]y^{-1}[y,w]\;.
$$
Clearly $[z,w]\in D_k(F;B\cup C\cup \{j\})= D_k(F;A\cup \{j\})$
and hence $y[z,w]y^{-1}\in D_k(F;A\cup \{j\})$. By inductive
hypothesis there exist $w'\in D_{k-1}(F;C\cup \{j\})$ and $x'\in
D_{k-1}(F;C)\cap F'$ such that $w=x'w'$. One has
$$
[y,w]=[y,x'w']=[y,x']x'[y,w']x^{\prime -1}\;.
$$
Clearly $[y,w']\in D_k(F;B\cup C\cup \{j\})= D_k(F;A\cup \{j\})$
and hence $x'[y,w']x^{\prime -1}\in D_k(F;A\cup \{j\})$. Therefore
there is an element $w''\in D_k(F;A\cup \{j\})$ such that
$[x,w]=[y,x']w''$ where $[y,x']\in D_k(F; A)\cap F'$.
\end{proof}

\medskip

For a given group $G$ the (lower) central series $\Gm_k=\Gm_k(G)$,
$$
G=\Gm_1\esq \Gm_2\esq \cdots \esq \Gm_k\esq \cdots
$$
of $G$ is defined inductively by
$\Gm_k=\uns{i+j=k}{\prod}[\Gm_i,\Gm_j]$. The well-known Witt-Hall
identities on commutators (see e.g. \cite{bc}) imply that
$\Gm_k=[G,\Gm_{k-1}]$.

Let ${\mathfrak Gr}$ denote the category of groups. Let us define
higher abelianization type functors $Z_k:{\mathfrak Gr}\to {\mathfrak
Gr}$, $k\ge 2$ by $Z_k(G)=G/\Gm_k(G)$ for any $G\in \ob{\mathfrak Gr}$
and where $Z_k(\al)$ is the natural homomorphism induced by a group
homomorphism $\al$.  Of course, $Z_2$ is the ordinary abelianization
functor of groups.

\medskip

%%Proposition 6
\begin{proposition}\label{2}
Let $(F; R_1, \ldots, R_n)$ be a normal $(n+1)$-ad of groups and
$k\ge 2$. Suppose that  $(F; R_1, \ldots, R_j)$ is a simple normal
$(j+1)$-ad of groups relative to $R_j$ for all $1\le j\le n$. Then
$$
\uns{i\in \la j\ra}{\cap}R_i \cap \Gm_k(F) = D_k(F; \la j\ra)\;,
\;\;\; 1\le j\le n\;.
$$
\end{proposition}
\begin{proof} Since the inclusion $D_k(F; \la j\ra)\seq
\uns{i\in \la j\ra}{\cap}R_i \cap \Gm_k(F)$ is clear, we only have
to show the inclusion $\uns{i\in \la j\ra}{\cap}R_i \cap
\Gm_k(F)\seq D_k(F; \la j\ra)$, which will be done by induction on
$j$.

Let $j=1$, then there exists a subgroup $F_1$ of the group $F$
such that $R_1\cap F_1=1$ and $F=F_1 R_1$. Let $w\in R_1\cap
\Gm_k(F)\seq \Gm_k(F)=D_k(F;\eset)$. Using Lemma \ref{1} one has
elements $x'\in D_k(F;\eset)\cap F_1$ and $w'\in D_k(F;\la 1\ra)$
such that $w=x' w'$. But $x'=ww^{'-1}\in R_1$ and hence $x'=1$.
Thus $R_1\cap \Gm_k(F)\seq D_k(F; \la 1 \ra)$.

Proceeding by induction, we suppose that the result is true for
$j-1$ and we will prove it for $j$.

There exists a subgroup $F_j$ of the group $F$ such that $R_j\cap
F_j = 1$ and $\uns{i\in A}{\cap} R_i = (\uns{i\in A}{\cap} R_i
\cap F_j)(\uns{i\in A}{\cap} R_i \cap R_j)$ for all $A\seq \{1,
\ldots, j-1\}$. Let $w\in \uns{i\in\la j\ra}{\cap}R_i\cap
\Gm_k(F)\seq \uns{i\in\la j-1 \ra}{\cap}R_i\cap \Gm_k(F)$. Using
the inductive hypothesis one has the equality $\uns{i\in \la j-1
\ra}{\cap}R_i\cap \Gm_k(F)=D_k(F;\la j-1\ra)$. By Lemma \ref{1}
there are elements $x'\in D_k(F;\la j-1\ra)\cap F_j$ and $w'\in
D_k(F;\la j\ra)$ such that $w=x'w'$. Certainly $x'=ww^{'-1}\in
R_j$ and hence $x'=1$. Therefore $\uns{i\in \la j\ra}{\cap}R_i\cap
\Gm_k(F)\seq D_k(F;\la j\ra)$.
\end{proof}

\medskip

These conditions may seem rather restrictive, but the following
observation shows that examples of simple normal $(n+1)$-ads of
groups appear naturally, and that moreover these examples satisfy
the conditions of Proposition \ref{2}.  First some terminology and
notation on pseudosimplicial groups, (cf. \cite{HI} for the
general theory).

A \emph{pseudosimplicial group} $G_*$ is a non-negatively graded group
with face   homomorphisms  $d^n_i:G_n\to G_{n-1}$ and
pseudodegeneracies, $s^n_i:G_n\to G_{n+1}$, $0\le i\le n$, satisfying
all the simplicial identities except possibly the identity $s^{n+1}_i
s^n_j=s^{n+1}_{j+1} s^n_i$ for $i\le j$ (again, see \cite{HI}). The Moore
complex $(NG_*,\pa_*)$ of $G_*$ is the chain complex defined by
$NG_n=\ovs{n-1}{\uns{i=0}{\cap}}\Ker d^n_i$ with $\pa_n:NG_n\to
NG_{n-1}$ induced from $d^n_n$ by restriction. The homotopy groups
of $G_*$ are defined as the homology groups of the complex
$(NG_*,\pa_*)$, i.e. $\pi_n(G_*)=H_n(NG_*,\pa_*)$, $n\ge 0$ (see
\cite{HI}). For good examples of pseudosimplicial groups see \cite{l2}.

\medskip

%%Proposition 7
\begin{proposition}\label{3}
Let $F_*$ be a pseudosimplicial group. Then $(F_n;\Ker
d^n_0,\ldots , \Ker d^n_{j-1})$ is a simple normal $(j+1)$-ad of
groups relative to $\Ker d^n_{j-1}$ for all $1\le j\le n$.
\end{proposition}
\begin{proof} Since $d^n_{j-1}s^{n-1}_{j-1}=1$,
$s^{n-1}_{j-1}(F_{n-1})\cap \Ker d^n_{j-1}=1$ and
$s^{n-1}_{j-1}(F_{n-1})\Ker d^n_{j-1}=F_n$ for all $n\ge 1$. Hence
when $j=1$, $(F_n;\Ker d^n_0)$ is a simple normal $2$-ad of groups
relative to $\Ker d^n_0$ and the $F^\prime$ of the definition of
simplicity is $s^{n-1}_0(F_{n-1})$.

Now suppose that $j>1$. We will show the following equality
$$
\uns{i\in A}{\cap}\Ker d^n_i =(\uns{i\in A}{\cap}\Ker d^n_i\cap
s^{n-1}_{j-1}(F_{n-1})) (\uns{i\in A}{\cap}\Ker d^n_i \cap \Ker
d^n_{j-1})
$$
for all $A\seq \{0,\ldots, j-2\}$ and $A\neq \eset$, so again the
$F^\prime$ of the definition of simplicity is
$s^{n-1}_{j-1}(F_{n-1})$.  Let $x=s^{n-1}_{j-1}(x_{n-1})r_{j-1}\in
\uns{i\in A}{\cap}\Ker d^n_i$, where $x_{n-1}\in F_{n-1}$,
$r_{j-1}\in \Ker d^n_{j-1}$. Thus $d^n_i(x)=d^n_i
s^{n-1}_{j-1}(x_{n-1})d^n_i(r_{j-1})=1$ for all $i\in A$. Since
$i< j-1$, one has
$d^n_i(r_{j-1})=s^{n-2}_{j-2}d^{n-1}_i(x_{n-1})^{-1}$. Hence
$1=d^{n-1}_id^n_{j-1}(r_{j-1})=d^{n-1}_{j-2}d^n_i(r_{j-1})=d^{n-1}_{j-2}
s^{n-2}_{j-2}d^{n-1}_i(x_{n-1})^{-1}=d^{n-1}_i(x_{n-1})^{-1}$.
Therefore $d^n_i(r_{j-1})=1$ and $d^n_i s^{n-1}_{j-1}(x_{n-1})=1$
for all $i\in A$.
\end{proof}

\medskip

The next lemma is well known, but very useful.  The proof is routine.

\medskip

%%Lemma 8
\begin{lemma}\label{4}
Let $G_*$ be a pseudosimplicial group and $A\seq \la n\ra$, $A\neq
\la n\ra$, then $d_n^n(\uns{i\in A}{\cap}\Ker d_{i-1}^n)=\uns{i\in
A}{\cap}\Ker d_{i-1}^{n-1}$, $n\ge 2$. \hfill$\Box$
\end{lemma}

\medskip

Next let us consider an augmented pseudosimplicial group $(F_*,
d^0_0, G)$ and, applying the functor $Z_k$ dimension-wise, denote
the resulting augmented pseudosimplicial group by $(Z_k(F_*),
Z_k(d^0_0), Z_k(G))$.  Our previous results will allow calculation
of the homotopy groups of $Z_k(F_*)$ in certain important cases
notably the following.

\medskip

%%Theorem 9
\begin{theorem}\label{5}
If $(F_*, d^0_0, G)$ is aspherical then there is a natural
isomorphism
$$
\pi_nZ_k(F_*)\cong \frac {\ovs{n-1}{\uns{i=0}{\cap}}\Ker
{d^{n-1}_i}\cap \Gm_k(F_{n-1})}{D_k (F_{n-1};\Ker d^{n-1}_0,
~\ldots,\Ker d^{n-1}_{n-1})}\;,\;\;\;k\ge 2,\;n\ge 1\;.
$$
\end{theorem}
\begin{proof} Let us consider the short exact sequence of
augmented pseudosimplicial groups
$$
\begin{matrix}
& 1& & 1& &1 & &1\\
& \da & & \da& &\da & &\da\\
\cdots \dtwo & \Gm_k(F_n) &\adtwo{\wt{d^n_0}}{\wt{d^n_n}} \cdots \three & \Gm_k(F_1) &
\atwo{\wt{d^1_0}}{\wt{d^1_1}} & \Gm_k(F_0) & \ovs{\wt{d_0^0}}{\to}& \Gm_k(G)\\
& \da & & \da& &\da & &\da\\
\cdots \dtwo & F_n &\adtwo{d^n_0}{d^n_n} \cdots \three & F_1 & \atwo{d^1_0}{d^1_1} & F_0 &
\ovs{d_0^0}{\to}& G\\
& \da & & \da& &\da & &\da\\
\cdots \dtwo & Z_k(F_n) &\dtwo \cdots \three & Z_k(F_1) & \two & Z_k(F_0) & \to & Z_k(G)\\
& \da & & \da& &\da & &\da\\
& 1& & 1& &1 & & 1
\end{matrix}
\;\;\;\;\;\;\;.
$$
By the induced long exact homotopy sequence, one has the
isomorphisms of groups
$$
\pi_nZ_k(F_*)\cong \frac{\ovs{n-1}{\uns{i=0}{\cap}}\Ker
\wt{d^{n-1}_i}}{
\wt{d^n_n}(\ovs{n-1}{\uns{i=0}{\cap}}\Ker\wt{d^n_i})},\;\quad n\ge
1\;.
$$

Since $\wt{d^n_i}$ is a restriction of $d^n_i$ to  $\Gm_k(F_n)$,
$\Ker\wt{d^n_i}=\Ker d^n_i\cap \Gm_k(F_n)$. Hence
$\ovs{n-1}{\uns{i=0}{\cap}}\Ker
\wt{d^{n-1}_i}=(\ovs{n-1}{\uns{i=0}{\cap}}\Ker {d^{n-1}_i})\cap
\Gm_k(F_{n-1})$ and $\ovs{n-1}{\uns{i=0}{\cap}}\Ker
\wt{d^n_i}=(\ovs{n-1}{\uns{i=0}{\cap}}\Ker{d^n_i})\cap
\Gm_k(F_n)$.

Using Proposition \ref{2} and Proposition \ref{3} one has
$$(\uns{i\in \la n\ra }{\cap}\Ker d^n_{i-1}) \cap \Gm_k(F_n)= \lbr
D_k (F_n;\Ker d^n_0,\ldots,\Ker d^n_{n-1})$$ for $n\ge 1$.

Since $(F_*, d^0_0, G)$ is an aspherical augmented
pseudosimplicial group, \lbr $d_n^n(\uns{i\in \la n \ra}{\cap}\Ker
d_{i-1}^n)=\uns{i\in \la n \ra}{\cap}\Ker d_{i-1}^{n-1}$, $n\ge
1$. Using this fact and Lemma \ref{4}, it is now easy to see that
one has an equality
$$
\wt{d^n_n}(\ovs{n-1}{\uns{i=0}{\cap}}\Ker \wt{d^n_i})=d^n_n(D_k
(F_n;\Ker d^n_0,\ldots,\Ker d^n_{n-1}))= \lbr D_k (F_{n-1};\Ker
d^{n-1}_0,\ldots,\Ker d^{n-1}_{n-1})\;.
$$
\end{proof}

\medskip

In the case where $(F_*, d^0_0, G)$ is a free pseudosimplicial
resolution of $G$, these homotopy groups will be the left
non-abelian derived functors $L_nZ_k(G)$ of $Z_k$, evaluated at
$G$.  (For more on these non-abelian derived functors, we refer
the reader to \cite{HI}.) We thus have the following formal
result.

\medskip

%%Corollary 10
\begin{corollary}\label{6}
Let $G$ be a group and $(F_*,d^0_0,G)$ an aspherical augmented
pseudosimplicial group and $k\ge 2$. If $F_n$ is a free group for
all $n\ge 0$, i.e. $(F_*,d^0_0,G)$ is a free pseudosimplicial
resolution of the group $G$, then there is a natural isomorphism
$$
L_nZ_k(G)\cong \frac {\ovs{n-1}{\uns{i=0}{\cap}}\Ker
{d^{n-1}_i}\cap \Gm_k(F_{n-1})}{D_k (F_{n-1};\Ker
d^{n-1}_0,\ldots,\Ker d^{n-1}_{n-1})}\;,\;\;\;n\ge 1\;.
$$
\end{corollary}
\begin{proof} Straightforward from Theorem \ref{5}.
\end{proof}

\

%%Section 3
\section{Crossed $n$-cubes of groups}

\

The following definition is due to Ellis and Steiner \cite{es}
(see also \cite{P}).  A \emph{crossed $n$-cube of groups} is a
family $\left\{{\mathcal M}_A:A\subseteq \la n \ra \right\}$ of
groups, together with homomorphisms, $\mu_i:{\mathcal M}_A\to
{\mathcal M}_{A\stmin \{i\}}$ for $i\in \la n \ra$, $A\seq \la n
\ra$ and functions $h:{\mathcal M}_A\times {\mathcal M}_B\lra
{\mathcal M}_{A\cup B}$ for $A,\;B\seq \la n \ra$, such that if
$^ab$ denotes $h(a,b)\cdot b$ for $a\in {\mathcal M}_A$ and $b\in
{\mathcal M}_B$ with $A\subseteq B$, then for all $a,\; a'\in
{\mathcal M}_A$, $b,\;b'\in {\mathcal M}_B$, $c\in {\mathcal M}_C$
and $i,\;j\in \la n \ra$, the following conditions hold:
\begin{align*}
& \mu_i(a)=a\;\;\text{if}\;\;i\notin A, \\
& \mu_i\mu_j(a)=\mu_j\mu_i(a),\\
& \mu_i h(a,b)=h(\mu_i(a),\mu_i(b)),\\
& h(a,b)=h(\mu_i(a),b)=h(a,\mu_i(b))\;\; \text{if}\;\; i\in A\cap B,\\
& h(a,a')=[a,a'],\\
& h(a,b)=h(b,a)^{-1},\\
& h(a,b)=1\;\; \text{if}\;\; a=1 \;\;\text{or}\;\; b=1,\\
& h(aa',b)= ~^ah(a',b)h(a,b),\\
& h(a,bb')=h(a,b) {}^bh(a,b'),\\
& {}^ah(h(a^{-1},b),c){}^ch(h(c^{-1},a),b){}^bh(h(b^{-1},c),a)=1,\\
& {}^ah(b,c)=h(^ab,^ac)\;\; \text{if}\;\; A\seq B\cap C.
\end{align*}

A morphism of crossed $n$-cubes, $\left\{{\mathcal M}_A
\right\}\to \left\{{\mathcal M}'_A \right\}$, is a family of group
homomorphisms $\left\{ f_A:{\mathcal M}_A\to {\mathcal M}'_A,\;
A\seq \la n \ra\right\}$ commuting with the $\mu_i$ and the
$h$-functions. This gives a category of crossed $n$-cubes of
groups which will be denoted by ${\it Crs}^n$.

\medskip

\

%%Examples
{\noindent \em Examples.}
\begin{enumerate}
\item[(i)] A crossed 1-cube is the same as a crossed module (see
Section 1). \item[(ii)] A crossed 2-cube is the same as a crossed
square (for the definition see \cite{bl}). The detailed
reformulation is easy. \item[(iii)] Let $G$ be a group and $N_1$,
$\ldots$ , $N_n$ be normal subgroups of $G$. Let ${\mathcal
M}_A=\uns{i\in A}{\cap} N_i$ for $A \seq \la n \ra$ (we also note
that then ${\mathcal M}_{\eset}=G$); if $i\in \la n\ra$, define
$\mu_i:{\mathcal M}_A \ovs{i\in A}{\lra}{\mathcal M}_{A\stmin
\{i\}}$ to be the inclusion and given $A,\;B\seq \la n \ra$, let
$h:{\mathcal M}_A\times {\mathcal M}_B\to {\mathcal M}_{A\cup B}$
be given by the commutator: $h(a,b)=[a,b]$ for $a\in {\mathcal
M}_A$, $b\in {\mathcal M}_B$  (here, of course,  ${\mathcal
M}_{A\cup B}={\mathcal M}_A\cap {\mathcal M}_B$) . Then
$\{{\mathcal M}_A: A\seq \la n \ra,\; \mu_i,\; h \}$ is a crossed
$n$-cube, called the \emph{inclusion crossed $n$-cube} given by
the \emph{normal $(n+1)$-ad of groups} $(G; N_1,\ldots, N_n)$.
\end{enumerate}

\medskip

Ellis and Steiner \cite{es} prove that ${\it Crs}^n$ is equivalent
to the category of $\cat ^n$-groups introduced by Loday who proved
that equivalence for $n=1, 2$,  \cite{l1}.

For a given crossed $n$-cube ${\mathcal M}$, there is an
associated ${\cat ^n}$-group and hence on applying the crossed
module nerve structure $E$, examined in Section 1, in the
$n$-independent directions, this construction leads naturally to
an $n$-simplicial group, called the \emph{multinerve} of the
crossed $n$-cube ${\mathcal M}$ and denoted by $\mathfrak
Ner({\mathcal M})$. Taking the diagonal of this $n$-simplicial
group gives a simplicial group denoted by $E^{(n)}({\mathcal
M})_*$, (see \cite{P}).

\medskip

%%Proposition 11
\begin{proposition}\label{7}
Let ${\mathcal M}$ be an inclusion crossed $n$-cube given by a
normal $(n+1)$-ad of groups $(F;R_1,\ldots,R_n)$ and $k\ge 2$.
Then there is a crossed $n$-cube ${\mathcal B}_k({\mathcal M})$
given by:
\begin{enumerate}
\item[(a)] for $A\seq \la n \ra$
$$
{\mathcal B}_k({\mathcal M})_A=\uns{i\in A}{\cap}R_i \diagup D_k(F;A)\;.
$$
\item[(b)] if $j\in \la n \ra$, the homomorphism
$\wt{\mu_j}:{\mathcal B}_k({\mathcal M})_A \to {\mathcal B}_k
({\mathcal M})_{A\stmin \{j\}}$ is induced by the inclusion
homomorphism $\mu_j$. \item[(c)] representing an element in
${\mathcal B}_k({\mathcal M})_A$ by $\ol x$ where $x\in \uns{i\in
A}{\cap}R_i$ (the bar denotes a coset), and for $A,\;B\seq \la n
\ra$, the map $\wt{h}:{\mathcal B}_k({\mathcal M})_A\times
{\mathcal B}_k({\mathcal M})_B\to {\mathcal B}_k({\mathcal
M})_{A\cup B}$ is given by $\wt h(\ol x,\ol
y)=\ol{h(x,y)}=\ol{[x,y]}$ for all $\ol x\in {\mathcal
B}_k({\mathcal M})_A$, $\ol y\in {\mathcal B}_k ({\mathcal M})_B$.
\end{enumerate}
\end{proposition}
\begin{proof} In our notation

$$ D_k(F;A)=\uns{A_1\cup A_2\cup \cdots \cup
A_k=A}{\prod}[\uns{i\in A_1}{\cap}R_i,[\uns{i\in
A_2}{\cap}R_i,\ldots,[\uns{i\in A_{k-1}}{\cap}R_i,\uns{i\in
A_k}{\cap}R_i]\ldots ]]\;, \quad A\seq \la n\ra. $$
Since \\
\begin{eqnarray*}[\uns{i\in
A_1}{\cap}R_i, [\uns{i\in A_2}{\cap}R_i,~ \ldots &&\hspace{-5mm},[\uns{i\in
A_{k-1}}{\cap}R_i,\uns{i\in A_k}{\cap}R_i]\ldots ]]\seq \\
&&[\uns{i\in A_1\stmin\{j\}}{\cap}R_i, [\uns{i\in A_2\stmin\{j\}}{\cap}R_i,
~\ldots ,[\uns{i\in A_{k-1}\stmin\{j\}}{\cap}R_i, \uns{i\in
A_k\stmin\{j\}}{\cap}R_i ]\ldots ]]
\end{eqnarray*}
 for $A_1\cup \cdots \cup
A_k=A\seq \la n\ra$, the inclusion $\mu_j:\uns{i\in A}{\cap}R_i
\hra \uns{i\in A\stmin\{j\}}{\cap}R_i$ induces the homomorphism
$\wt{\mu_j}:{\mathcal B}_k({\mathcal M})_A\to {\mathcal B}_k({\mathcal
M})_{A\stmin\{j\}}$ for all $j\in \la n\ra$.

Now we are only left to show that the function $\wt{h}:{\mathcal
B}_k({\mathcal M})_A\times {\mathcal B}_k({\mathcal M})_B\to
{\mathcal B}_k({\mathcal M})_{A\cup B}$ for $A,B\seq \la n\ra$ is
well defined. In fact, let $x'\in \uns{i\in A}{\cap}R_i$,
$y'\in \uns{i\in B}{\cap}R_i$ be  such that $$xx^{\prime -1}\in
\uns{A_1\cup \cdots \cup A_k=A}{\prod}[\uns{i\in A_1}{\cap}R_i,
[\uns{i\in A_2}{\cap}R_i, \ldots ,[\uns{i\in
A_{k-1}}{\cap}R_i,\uns{i\in A_k}{\cap}R_i]\ldots ]]$$ and
$$yy^{\prime -1}\in \lbr \uns{A_1\cup \cdots \cup
A_k=B}{\prod}[\uns{i\in A_1}{\cap}R_i, [\uns{i\in A_2}{\cap}R_i,
\ldots ,[\uns{i\in A_{k-1}}{\cap}R_i,\uns{i\in
A_k}{\cap}R_i]\ldots ]].$$ The inclusion $$[\uns{i\in A}{\cap}R_i,
\uns{i\in B}{\cap}R_i]\seq \uns{i\in A\cup B}{\cap}R_i$$ for all
$A,B\seq \la n\ra$ implies that

\begin{eqnarray*}
[x,y][x',y']^{-1}&=& xyx^{-1}y^{-1}y'x'y^{\prime -1}x^{\prime -1}\\
&=& xy'[y^{\prime -1}y,x^{-1}] y^{\prime -1}x^{-1}x
[y',x^{-1}x']x^{-1}\\ &\in& \uns{A_1\cup \cdots \cup A_k=A\cup
B}{\prod}[\uns{i\in A_1}{\cap}R_i, [\uns{i\in A_2}{\cap}R_i,
\ldots , [\uns{i\in A_{k-1}}{\cap}R_i,\uns{i\in
A_k}{\cap}R_i]\ldots ]]\;,
\end{eqnarray*}
so $\ol{h(x,y)}=\ol{h(x',y')}$ and $\wt h$ is well defined. The
verification that $\mathcal{B}_k(\mathcal{M})$ is a crossed
$n$-cube is routine and is left as an exercise.
\end{proof}

\medskip

%%Remark
{\noindent \em Remark.} The functor ${\mathcal B}_2$ coincides on
the subcategory of inclusion crossed $n$-cubes with the
abelianization functor ${\it Ab}$ from the category ${\it Crs}^n$
to the category ${\mathfrak Ab}{\it Crs}^n$ of abelian crossed
$n$-cubes of groups (i.e. crossed $n$-cubes all of whose $h$ maps
are trivial), considered for $n=1$ in Section 1 and defined in
general by the following way: for any ${\mathcal M}=\{{\mathcal
M}_A:A\seq \la n \ra,\; \mu_i,\; h\}$ in $ {\it Crs}^n$, ${\it
Ab}(\mathcal M)$ is an abelian crossed $n$-cube given by:
\begin{enumerate}
\item[(a)] for $A\seq \la n \ra$
$$
{\it Ab}({\mathcal M})_A=\frac{{\mathcal M}_A}{ \uns{B\cup
C=A}{\uns{B,C\; \seq \la n\ra}{\prod}} D_{B,C}}\;,
$$
where $D_{B,C}$ is the subgroup of ${\mathcal M}_A$ generated by
the elements $h(b,c)$, \lbr $h:{\mathcal M}_B\times {\mathcal M}_C
\to {\mathcal M}_{B\cup C=A}$ for all $b\in {\mathcal M}_B$, $c\in
{\mathcal M}_C$. \item[(b)] if $i\in \la n \ra$, the homomorphism
$\wt\mu_i:{\it Ab}({\mathcal M})_A \to {\it Ab} ({\mathcal
M})_{A\stmin \{ i \}}$ is induced by the homomorphism $\mu_i$.
\item[(c)] for $A,\;B\seq \la n \ra$, the function $\wt{h}:{\it
Ab}({\mathcal M})_A \times {\it Ab}({\mathcal M})_B\to {\it
Ab}({\mathcal M})_{A\cup B}$ is induced by $h$ and therefore is
trivial.
\end{enumerate}
The functor ${\it Ab}:{\it Crs}^n \to {\mathfrak Ab}{\it Crs}^n$ is
left adjoint to the inclusion functor \lbr $i:{\mathfrak Ab}{\it
Crs}^n\hra {\it Crs}^n$ as is easily checked.

\medskip

For any inclusion crossed $n$-cube ${\mathcal M}$ given by a
normal $(n+1)$-ad of groups $(F;R_1,\ldots,R_n)$ and $k\ge 2$,
there is a natural morphism of crossed $n$-cubes ${\mathcal M}\to
{\mathcal B}_k({\mathcal M})$ inducing the natural fibration of
simplicial groups $E^{(n)}({\mathcal M})_*\ovs{\De^{n,k}_*}{\to}
E^{(n)}({\mathcal B}_k({\mathcal M}))_*$ defined by
$$\De^{n,k}_m(x_1,\ldots, x_l)=(\ol{x_1},\ldots, \ol{x_l})$$ for all
$(x_1,\ldots, x_l)\in E^{(n)}({\mathcal M})_m = (\uns{i\in
A_1}{\cap} R_i)\rtimes \cdots \rtimes (\uns{i\in A_l}{\cap} R_i)$
and $m\ge 0$, where $A_1,\ldots,A_l\seq \la n \ra$ and
$l={(m+1)}^n$. It is easy to see that $\Ker
\De^{n,k}_m=D_k(F;A_1)\rtimes \cdots \rtimes D_k(F;A_l)$.

Let us consider a simplicial group $G_*$ and, applying the functor
$Z_k$ dimension-wise, denote the resulting  simplicial group by
$Z_k G_*$.

\medskip

%%Proposition 12
\begin{proposition}\label{8}
Let ${\mathcal M}$ be an inclusion crossed $n$-cube given by a
normal $(n+1)$-ad of groups $(F;R_1,\ldots, R_n)$  and $k\ge 2$.
Then there is an isomorphism of simplicial groups
$$
Z_kE^{(n)}({\mathcal M})_*\cong E^{(n)}({\mathcal B}_k({\mathcal M}))_*\;.
$$
\end{proposition}
\begin{proof}
For any inclusion crossed module $R\hra F$, It is easy to check
the following equalities in the group $R\rtimes \cdots \rtimes
R\rtimes F$:
\begin{align*}
&[(1,\ldots,1,x),(1,\ldots,1,x')]=(1,\ldots,1,[x,x'])\;,\\
&[(1,\ldots,1,\uns{s}{r},1,\ldots,1),(1,\ldots,1,x)]=(1,\ldots,1,\uns{s}{[r,x]},1,\ldots,1)\;,\\
&[(1,\ldots,1,\uns{s}{r},1,\ldots,1),(1,\ldots,1,\uns{t}{r'},1,\ldots,1)]=
(1,\ldots,1,\uns{\min\{s,\;t\}}{[r,r']},1,\ldots,1)
\end{align*}
for all $x,x'\in F$, $r,r'\in R$.

There are further generalisation of these equalities, namely for
any inclusion crossed $n$-cube $\mathcal M$ given by the normal
$n+1$-ad of groups $(F,R_1,\ldots, R_n)$ one has the following
facts, the proof of which is routine and will be omitted.

\medskip

\begin{enumerate}
\item[(A)] Let $s$ and $t$ be any fixed elements of the set $\la
(m+1)^n\ra$. Then there exists a unique $\lam=\lam(s, t)\in \la
(m+1)^n\ra$ such that $A_{\lam}=A_s\cup A_t$ and in the group
$E^{(n)}({\mathcal M})_m$ there holds the equality:
$$
[(1,\ldots,1,\uns s x, 1, \ldots, 1),(1,\ldots,1,\uns t y, 1,
\ldots, 1)]=(1,\ldots,1,\uns \lam {[x,y]}, 1, \ldots, 1)
$$
for all $x\in \uns{i\in A_s}{\cap} R_i$, $y\in \uns{i\in
A_t}{\cap} R_i$.

\item[(B)] Let $s\in \la (m+1)^n\ra$ and $A, B\seq A_s$ with
$A\cup B=A_s$. Then there exists $p, q \in \la (m+1)^n\ra$ such
that $A_p=A$, $A_q=B$ and $\lam(p,q)=s$.
\end{enumerate}

\medskip

We only have to show the equality
\begin{equation}\label{X}
\Gm_k(E^{(n)}({\mathcal M})_m) = \Ker \De^{n,k}_m
\end{equation}
which will be done by induction on $k$, using facts (A) and (B)
above.

Let $k=1$, then it is clear that $\Gm_1(E^{(n)}({\mathcal M})_m) =
\Ker \De^{n,1}_m$.

Proceeding by induction, we suppose that (\ref{X}) is true for
$k-1$ and we will prove it for $k$.

First we will show the inclusion $\Ker \De^{n,k}_m\seq
\Gm_k(E^{(n)}({\mathcal M})_m)$. It is sufficient to show
$$
1\rtimes \cdots \rtimes 1\rtimes D_k(F,A_s)\rtimes 1\rtimes \cdots \rtimes 1\seq
\Gm_k(E^{(n)}({\mathcal M})_m) \quad \text{for all}\quad s\in \la
(m+1)^n\ra\;.
$$
In fact, any generator $w$ of $D_k(F,A_s)$ has the form $w=[x,y]$,
where $x\in \uns{i\in A}{\cap}R_i$, $y\in D_{k-1}(F,B)$ and $A\cup
B=A_s$.

Now (B) implies that there exist $p,q\in \la (m+1)^n\ra$ such that
$A_p=A$, $A_q=B$ and $\lam(p,q)=s$. Thus we have
\begin{align*}
&[(1,\ldots,1,\uns p x, 1, \ldots, 1), (1,\ldots,1,\uns q y, 1,
\ldots, 1)]=(1,\ldots,1,\uns s w, 1,\ldots,1)\;,
\end{align*}
which means that
$$
1\rtimes \cdots \rtimes 1\rtimes D_k(F, A_s)\rtimes 1\rtimes
\cdots \rtimes 1\seq [E^{(n)}({\mathcal M})_m, \Ker
\De^{n,k-1}_m]\;.
$$
Therefore by the inductive hypothesis we obtain
$$
1\rtimes \cdots \rtimes 1\rtimes D_k(F, A_s)\rtimes 1\rtimes
\cdots \rtimes 1\seq [E^{(n)}({\mathcal M})_m,
\Gm_{k-1}(E^{(n)}({\mathcal M})_m)]= \Gm_k(E^{(n)}({\mathcal
M})_m)\;.
$$

Finally we will show the inverse inclusion
$\Gm_k(E^{(n)}({\mathcal M})_m)\seq \Ker \De^{n,k}_m$. In fact,
any generator $w$ of $\Gm_k(E^{(n)}({\mathcal M})_m)$ could be
written as $w=[w_1,w_2]$, where $w_1\in E^{(n)}({\mathcal M})_m$
and $w_2\in \Gm_{k-1}(E^{(n)}({\mathcal M})_m)$. Using again the
inductive hypothesis we have $w_2\in \Ker \De^{n,k-1}_m$. Thus
\begin{align*}
&w_1=\ovs{(m+1)^n}{\uns{s=1}{\prod}}(1,\ldots, 1, \uns{s}{x_s}, 1,
\ldots, 1)\;,\quad x_s\in \uns{i\in A_s}{\cap} R_i\;,\\
&w_2=\ovs{(m+1)^n}{\uns{t=1}{\prod}}(1,\ldots, 1, \uns{t}{y_t}, 1,
\ldots, 1)\;,\quad y_t\in D_{k-1}(F, A_t)\;.
\end{align*}
It is certain that $[x_s, y_t]\in D_k(F, A_s\cup A_t)$. Then (A)
implies that we have
$$
\begin{matrix}
[(1,\ldots,1,\uns{s}{x_s}, 1, \ldots, 1),
(1,\ldots,1,\uns{t}{y_t}, 1, \ldots,
1)]=(1,\ldots,1,\uns{\lam(s,t)}{[x_s,y_t]}, 1,\ldots,1)\\ \in
1\rtimes \cdots \rtimes 1\rtimes D_k(F, A_{\lam(s,t)})\rtimes
1\rtimes \cdots \rtimes 1\seq \Ker \De^{n,k}_m\;.
\end{matrix}
$$
Then the Witt-Hall identities on commutators implies that $w\in
\Ker \De^{n,k}_m$.
\end{proof}

\

%%section 4
\section{Non-abelian mapping cone complex}

\

A complex of (non-abelian) groups $(A_*,d_*)$ of length $n$ is a
sequence of group homomorphisms
$$
A_n \ovs{d_n}{\lra} A_{n-1}\ovs{d_{n-1}}{\lra} \cdots \ovs{d_1}{\lra} A_0
$$
such that $\Im d_{i+1}$ is normal in $\Ker d_i$. Now we recall the
following definition from \cite{l1}.

Let $f:(A_*,d_*)\to (B_*,d'_*)$ be a morphism of chain complexes
of groups. Let $f$ satisfy the following conditions  ($\ast$):

\emph{each} $f_i: A_i \to B_i$ \emph{is a crossed module} \\and

\emph{the maps $(d_i, d'_i)$
form a morphism of crossed modules.}

Then the mapping cone of $f$ is a complex of (non-abelian) groups
$(C_*(f),\pa _*)$ defined by $C_i(f)=A_{i-1}\rtimes B_i$, where
the action of $B_i$ on $A_{i-1}$ is induced by the action of
$B_{i-1}$ on $A_{i-1}$ via the homomorphism $d'_i$; and
$$\pa_i(a,b)=(d_{i-1}(a)^{-1},f_{i-1}(a)d'_i(b))$$ for all $a\in
A_{i-1}$, $b\in B_i$. Then by \cite{l1}, Proposition 3.2, there is
a long exact sequence of groups.
\begin{equation}\label{A}
\cdots \lra H_i(A_*)\lra H_i(B_*)\lra H_i(C_*(f))\lra
H_{i-1}(A_*)\lra \cdots \;.
\end{equation}

Now let us consider a morphism of pseudosimplicial groups
$\al:(G_*,d^*_i,s^*_i)\to (H_*,d^{\prime*}_i,s^{\prime*}_i)$
satisfying conditions ($\ast$ $\ast$)  :

\emph{each} $\al_n:G_n \to H_n$ \emph{is a crossed module} \\and

\emph{the maps $(d^*_i, d^{\prime *}_i)$ and $(s^*_i,s^{\prime
*}_i)$ form  morphisms of crossed modules.}

Define a new pseudosimplicial group $M_*(\al)$ in the following
way:
$$
M_n(\al)=\underbrace{G_n\rtimes G_n \rtimes \cdots \rtimes
G_n}_{n-\text{times}} \rtimes H_n\;,
$$
\begin{align*}
&d^n_0(g_1,\ldots,g_n,h)=(d^n_0(g_2),\ldots, d^n_0(g_n),d^{\prime n}_0(h))\;,\\
&d^n_i(g_1,\ldots,g_n,h)=(d^n_i(g_1),\ldots, d^n_i(g_i) d^n_i(g_{i+1}),\ldots,d^n_i(g_n),
d^{\prime n}_i(h))\;,\;\;\;0<i<n\;,\\
&d^n_n(g_1,\ldots,g_n,h)=(d^n_n(g_1),\ldots,
d^n_n(g_{n-1}),\al_{n-1}d^n_n(g_n)
d^{\prime n}_n(h))\;,\\
&s^n_i(g_1,\ldots,g_n,h)=(s^n_i(g_1), \ldots, s^n_i(g_i), 1,
s^n_i(g_{i+1}),\ldots, s^n_i(g_n), s^{\prime n}_i(h))\;,\;\;\;0\le
i\le n\;.
\end{align*}

It is easy to see that the induced morphism $\wt\al:NG_*\to NH_*$,
where $NG_*$ and $NH_*$ are the Moore complexes of $G_*$ and $H_*$
respectively, satisfies the conditions ($\ast$). Therefore one can
consider the mapping cone complex $C_*(\wt\al)$ of $\wt\al$.

\

%%Proposition 13
\begin{proposition}\label{9}
The natural morphism of complexes $\ka: NM_*(\al)\to C_*(\wt\al)$,
given by $\ka_n(g_1,g_2,\ldots,g_n,h)=(d^n_n(g_n),h)$, $n\ge 0$
induces an isomorphism of groups
$$
\pi_n(M_*(\al))\cong H_n(C_*(\wt\al))\;,\quad n\ge 0\;.
$$
\end{proposition}
\begin{proof} The verification that $\ka_n$, $n\ge 0$
is a homomorphism and commuting with differentials is easy. Let
$(g,h)\in NG_{n-1}\rtimes NH_n=C_n(\wt\al)$, then it is easy to
check that
$(s^{n-1}_0(g)^{\ep(n-1)},\ldots,s^{n-1}_{n-2}(g)^{-1},s^{n-1}_{n-1}(g),h)\in
NM_n(\al)$, where $\ep(i)=(-1)^i$. It is clear that
$\ka_n(s^{n-1}_0(g)^{\ep(n-1)},\ldots,s^{n-1}_{n-2}(g)^{-1},s^{n-1}_{n-1}(g),h)=(g,h)$.
Hence $\ka_n$ is surjective for all $n\ge 0$.

Consider the kernel complex $({\mathfrak G}_*, \ol{\pa}_*)$ of $\ka$.
Note that $\Im \ol{\pa}_n$ is not normal in $\Ker \ol{\pa}_{n-1}$
in general, ${\mathfrak G}_0=1$ and
$$
{\mathfrak G}_n=\left\{\begin{matrix}(g_1,g_2,\ldots,g_n)\in\\
\underbrace{G_n\rtimes G_n \rtimes \cdots \rtimes
G_n}_{n-\text{times}} \end{matrix}
\begin{matrix}
|&d^n_0(g_j)=1\;, \quad 2\le j\le n\; ;\\
|&d^n_i(g_j)=d^n_i(g_i)d^n_i(g_{i+1})=1\;,\quad 1\le i\le
n-1\;,\\ |&\hfill 1\le j\le n\;\;\text{and}\;\; i\neq j-1, j\;;\\
|&d^n_n(g_n)=1\;
\end{matrix}
\right\}\;\;.
$$
Furthermore, it is easy to check that for an element $(g_1,~~\ldots,~~
g_{n-1})\in \Ker \ol{\pa_{n-1}}$, the element $(g'_1, \ldots,
g'_{n-1}, g'_n)$, defined by the formulae:
$$
g'_i = \left\lbrace
  \begin{array}{cl}
s^{n-1}_{n-1}(g_i) s^{n-1}_{n-2}(g^{-1}_i)\cdots s^{n-1}_i(g^{\ep(n-i-1)}_i) s^{n-1}_{i-1}(g^{\ep(n-i)}_i g^{\ep(n-i)}_{i+1}\cdots g^{\ep(n-i)}_{n-1}),\; i\;\text{even},\\
s^{n-1}_{i-1}( g^{\ep(n-i)}_{n-1}\cdots g^{\ep(n-i)}_{i+1} g^{\ep(n-i)}_i) s^{n-1}_i(g^{\ep(n-i-1)}_i)\cdots s^{n-1}_{n-2}(g^{-1}_i) s^{n-1}_{n-1}(g_i),\; i\;\text{odd},\\
 \end{array} \right.
$$
for all $1\le i\le n-1$ and $g'_n=1$, belongs to ${\mathfrak G}_n$ and
$$\ol{\pa_n}(g'_1, \ldots, g'_{n-1}, g'_n)=(g_1, \ldots, g_{n-1}).$$
Now the proposition  follows from the long exact homology sequence
induced by the short exact sequence of complexes $1\lra {\mathfrak
G}_*\lra NM_*(\al)\ovs{\ka}{\lra}C_*(\wt\al)\lra 1$.
\end{proof}

\medskip

Given a pseudosimplicial group $G_*$, we will say that the \emph{length
of $G_*$ is }$\le n$, denoted by $l(G_*)\le n$ if $NG_i=1$ for $i>n$.

\medskip

%%Remark
{\noindent \em Remark.} Let $\al:(G_*,d^*_i,s^*_i)\to
(H_*,d^{\prime*}_i,s^{\prime*}_i)$ be a morphism of
pseudosimplicial groups satisfying the conditions ($\ast$$\ast$)
and $n\ge 2$. Suppose $l(G_*)\le n-1$ and $l(H_*)\le n-1$. Consider an
element $(g_1,g_2,\ldots,g_k,h)\in NM_k(\al)$, $k>n$, then
\begin{align*}
&d^k_0(g_j)=1\;,\quad 2\le j\le k\;,\\
&d^k_i(g_j)=d^k_i(g_i)d^k_i(g_{i+1})=1\;,\quad 1\le i\le
k-1\;,\;\;1\le j\le k\;\;\text{and}\;\; i\neq j-1\;, j\;,\\
&d^{\prime k}_i(h)=1\;,\quad 0\le i\le k-1.
\end{align*}
By Lemma \ref{4} one can easily show that $g_i=1$, $1\le i\le k$
and $h=1$, meaning  $NM_k(\al)=1$ for $k>n$. Thus $l(M_*(\al))\le
n$.

\medskip

Now using the mapping cone construction, for a given crossed
$n$-cube ${\mathcal M}$, we construct inductively a complex of
groups $C_*({\mathcal M})$ of length $n$, always having in mind
that ${\mathcal M}$ is thought as a crossed module of crossed
$(n-1)$-cubes, ${\mathcal M}_1\to {\mathcal M}_0$ (Proposition 5,
\cite{P}). In fact, for $n=1$, and ${\mathcal
M}=(M\ovs{\mu}{\to}P)$, $C_*(\mathcal{M})$ is the complex $M\to P$
of length $1$. Let $n=2$ and ${\mathcal M}$ be a crossed square,
considered as a crossed module of crossed modules or a morphism of
complexes of length $1$ satisfying the conditions ($\ast$). The
construction above gives a complex $C_*({\mathcal M})$ of length
$2$. (This has a 2-crossed module structure, \cite{dc}, as noted
by Conduch\'{e}, see also \cite{atmp5}.) Proceeding by induction,
suppose for any crossed $(n-1)$-cube ${\mathcal M}$ we have
constructed a complex $C_*({\mathcal M})$ of length $n-1$. Now let
${\mathcal M}$ be a crossed $n$-cube and consider  it as a crossed
module of crossed $(n-1)$-cubes ${\mathcal M}_1\to {\mathcal
M}_0$, which implies a morphism of complexes of groups
$C_*({\mathcal M}_1)\ovs{\de}{\lra} C_*({\mathcal M}_0)$ of length
$n-1$ satisfying the conditions ($\ast$). So using again the
above-mentioned construction we obtain a complex of groups
$C_*({\mathcal M})=C_*(\de)$ of length $n$.

\medskip

%%Proposition 14
\begin{proposition}\label{10}\cite{l1}
Let ${\mathcal M}$ be a crossed $n$-cube of groups. Then
$l(E^n({\mathcal M})_*)\le n$ and there is a natural morphism of
complexes $NE^{(n)}({\mathcal M})_*\to C_*({\mathcal M})$ which
induces isomorphisms of groups
$$
\pi_i(E^{(n)}({\mathcal M})_*)\cong H_i(C_*({\mathcal M}))\;,\quad
i\ge 0\;.
$$
Moreover $\pi_n(E^{(n)}({\mathcal M})_*)\cong
\ovs{n}{\uns{i=1}{\cap}} \Ker({\mathcal M}_{\la n
\ra}\ovs{\mu_i}{\lra} {\mathcal M}_{\la n \ra \ \stmin \{i\}})$.
\end{proposition}
\begin{proof} This is obvious for $n=1$. Let $n=2$ and
${\mathcal M}$ be a crossed square respectively. If we consider
${\mathcal M}$ as a crossed module of crossed modules ${\mathcal
M}_1\to {\mathcal M}_0$ inducing the natural morphism of
simplicial groups $E^{(1)}({\mathcal M}_1)_*\ovs{\al}{\to}
E^{(1)}({\mathcal M}_0)_*$ satisfying the conditions
($\ast$$\ast$), then by definition $E^{(2)}({\mathcal
M})_*=M_*(\al)$, and by Proposition \ref{9} and the corresponding
{\em Remark}, $l(E^{(2)}({\mathcal M})_*)\le 2$, and there exists
a natural morphism of complexes $NE^{(2)}({\mathcal M})_*\to
C_*(\wt\al)$ inducing an isomorphism $$\pi_i(E^{(2)}({\mathcal
M})_*)\cong H_i(C_*(\wt\al))\;,\quad i\ge 0\;.$$ Clearly
$C_*(\wt\al)\cong C_*({\mathcal M})$.

Proceeding by induction, we suppose that the assertion is true for
$n-1$ and we will show it for $n$.

Let us consider any crossed $n$-cube ${\mathcal M}$ as a crossed
module of crossed $(n-1)$-cubes ${\mathcal M}_1\to {\mathcal
M}_0$. This implies a morphism of simplicial groups
$E^{(n-1)}({\mathcal M}_1)_*\ovs{\al}{\to} E^{(n-1)}({\mathcal
M}_0)_*$ satisfying the conditions ($\ast$$\ast$) and a morphism
of complexes $C_*({\mathcal M}_1)\ovs{\de}{\to} C_*({\mathcal
M}_0)$ satisfying the conditions ($\ast$). By definition
$E^{(n)}({\mathcal M})_*=M_*(\al)$, hence Proposition \ref{9} and
its {\em Remark}  imply that $l(E^{(n)}({\mathcal M})_*)\le n$ and
there exists a natural morphism of complexes $NE^{(n)}({\mathcal
M})_*\ovs{\ka}{\to} C_*(\wt\al)$ inducing isomorphisms
$$\pi_i(E^{(n)}({\mathcal M})_*)\cong H_i(C_*(\wt\al)), \quad i\ge
0.$$ Using the inductive hypothesis, there exist  natural morphisms
of complexes\\
\centerline{ $NE^{(n-1)}({\mathcal M}_1)_*\ovs{\ka'}{\to}
C_*({\mathcal M}_1)$ \quad and \quad $NE^{(n-1)}({\mathcal
M}_0)_*\ovs{\ka''}{\to} C_*({\mathcal M}_0)$,}
 which induce
isomorphisms $$\pi_i(E^{(n-1)}({\mathcal M}_1)_*)\cong
H_i(C_*({\mathcal M}_1)),$$  $$\pi_i(E^{(n-1)}({\mathcal
M}_0)_*)\cong H_i(C_*({\mathcal M}_0)),$$ for  $i\ge 0$. It is
easy to check that $\ka''\wt\al=\de\ka'$ and that
$(\ka'_i,\ka''_i)$ is a morphism of crossed modules for all $i\ge
0$. Then the natural morphism of complexes
$C_*(\wt\al)\ovs{\ka'\rtimes\ka''}{\lra}C_*(\de)=C_*({\mathcal
M})$, by (\ref{A}) and the  five lemma, induces
$H_i(C_*(\wt\al))\cong H_i(C_*({\mathcal M}))$, $i\ge 0$.
Therefore the morphism of complexes $$NE^{(n)}({\mathcal
M})_*\ovs{(\ka'\rtimes\ka'')\circ \ka}{\lra}C_*({\mathcal M})$$
induces $$\pi_i(E^{(n)}({\mathcal M})_*)\cong H_i(C_*({\mathcal
M})),  \quad i\ge 0.$$ From these isomorphisms and the
construction of $C_*({\mathcal M})$ follows that
$$\pi_n(E^{(n)}({\mathcal M})_*)\cong \ovs{n}{\uns{i=1}{\cap}}
\Ker({\mathcal M}_{\la n \ra}\ovs{\mu_i}{\lra} {\mathcal M}_{\la n
\ra \ \stmin \{i\}})\;.$$
\end{proof}

\

%%Section 5
\section{$n$-fold \v{C}ech derived functors}

\

The \v{C}ech derived functors of group valued functors were
introduced in \cite{Pi1} (see also \cite{HI} and, here, our
Section 1) as an algebraic analogue of the \v{C}ech (co)homology
construction of open covers of topological spaces with
coefficients in sheaves of abelian groups (see \cite{Sh}). It is
well known that the \v{C}ech cohomlogy of topological spaces with
coefficients in sheaves is closely related to sheaf cohomology of
topological spaces, in particular this relation is expressed by
spectral sequences \cite{Sh}. Some applications of \v{C}ech
derived functors to group (co)homology theory and $K$-theory are
given in \cite{Pi1, Pi2, Pi3}. In this section we generalise the
notion of the \v{C}ech derived functors to that of the $n$-fold
\v{C}ech derived functors of an endofunctor on the category of
groups. Based on this notion we get a new purely algebraic method
for the investigation of higher integral homology of groups from a
Hopf formula point of view and the further generalizations of
these formulae.

Let us consider again the set $\la n \ra = \{1, \ldots, n\}$. The
subsets of $\la n \ra$ are ordered by inclusion. This ordered set
determines in the usual way a category $\ul{C_n}$. For every pair
$(A,B)$ of subsets with $A\seq B\seq \la n \ra$, there is the
unique morphism $\rho^A_B:A\to B$ in $\ul{C_n}$. It is easy to see
that any morphism in the category $\ul{C_n}$, not an identity, is
generated by $\rho^A_B$ for all $A\seq \la n \ra$, $A\neq \la n
\ra$, $B = A\cup \{j\}$, $j\notin A$.

An \emph{$n$-cube of groups} is a functor ${\mathfrak
F}:\ul{C_n}\to {\mathfrak Gr}$, $A\mapsto {\mathfrak F}_A$,
$\rho^A_B\mapsto \al^A_B$. A \emph{morphism between $n$-cubes}
${\mathfrak F}, {\mathfrak Q}:\ul{C_n}\to {\mathfrak Gr}$ is a
natural transformation $\ka:{\mathfrak F}\to {\mathfrak Q}$.

\medskip

{\sc Warning:} A crossed $n$-cube of groups gives an $n$-cube on
forgetting structure, but note that there is a reversal of the
role of the index $A$.  The top corner of a crossed $n$-cube is
$G_{\la n \ra}$, that in an $n$-cube is $\mathfrak{F}_\emptyset$.
This is due to the fact that an $n$-cube of groups yields a
crossed $n$-cube as a sort of generalized kernel.

\medskip

Let $A\seq \la n\ra$ and consider two full subcategories of the
category $\ul{C_n}$: $\ul{C_n^A}$ is the category of all subsets
of $\la n\ra$ containing the subset $A$ and $\ul{C_n^{\ol{A}}}$ is
the category of all subsets of $\la n\ra$ not containing the
subset $A$. For a given $n$-cube of groups ${\mathfrak F}$, and
$A$ as above, denote by ${\mathfrak F}^A$ and ${\mathfrak
F}^{\ol{A}}$ the functors induced by the restriction of the
functor ${\mathfrak F}$ to the subcategories $\ul{C_n^A}$ and
$\ul{C_n^{\ol{A}}}$ respectively. For a given morphism of
$n$-cubes of groups $\ka:{\mathfrak F}\to {\mathfrak Q}$ denote by
$\ka^A:{\mathfrak F}^A \to {\mathfrak Q}^A$ the natural
transformation induced by restriction of the natural
transformation $\ka$.

\medskip

%%Examples
{\noindent \em Examples}
\begin{enumerate}
\item[(a)] Let $(F;R_1,\ldots, R_n)$ be a normal $(n+1)$-ad of
groups. These data naturally determines an $n$-cube of groups
${\mathfrak F}$ as follows: for any $A\seq \la n \ra$, let
${\mathfrak F}_A=F\diagup \uns{i\in A}{\prod}R_i$; for an
inclusion $A\seq B$, let $\al^A_B:{\mathfrak F}_A\to {\mathfrak
F}_B$ be the natural homomorphism induced by $1_F$. This $n$-cube
of groups will be called the \emph{$n$-cube of groups induced by
the normal $(n+1)$-ad of groups}, $(F;R_1,\ldots, R_n)$.
\item[(b)] Let $(G_*, d^0_0, G)$ be an augmented pseudosimplicial
group. A natural $n$-cube of groups ${\mathcal G}(n):\ul{C_n}\to
{\mathfrak Gr}$, $n\ge 1$ is defined by the following way:
$$
\begin{matrix}
{\mathcal G}(n)_A=G_{n-1-|A|}\;\;\text{for all}\; A\seq \la n\ra \;,\\
\al^A_{A\cup \{j\}}= d^{n-1-|A|}_{k-1}\;\;\text{for all}\; A\neq
\la n\ra\;,\;j\notin A\;,
\end{matrix}
$$
where $G_{-1}=G$, $\de(k)=j$ and $\de: \la n-|A|\ra \to \la n\ra
\stmin A$ is the unique monotone bijection.
\end{enumerate}

\medskip

Given an $n$-cube of groups ${\mathfrak F}$. It is easy to see
that there exists a natural homomorphism ${\mathfrak F}_A
\ovs{\al_A}{\lra}\uns {B\supset A}{\lim}{\mathfrak F}_B$ for any
$A\seq \la n \ra$, $A\neq \la n\ra$.

Let $G$ be a group. An $n$-cube of groups ${\mathfrak F}$ will be
called an \emph{$n$-presentation of the group} $G$ if ${\mathfrak
F}_{\la n\ra}=G$. An $n$-presentation ${\mathfrak F}$ of $G$ is
called \emph{free} if the group ${\mathfrak F}_A$ is free for all
$A\neq \la n \ra$ and called \emph{exact} if the homomorphism
$\al_A$ is surjective for all $A\neq \la n\ra$. Note that a
fibrant $n$-presentation of a group $G$ in the sense of
Brown-Ellis \cite{be} is the same as a free exact $n$-presentation
of $G$ in our sense, for a construction of such, see \cite{be}.

\medskip

%%Proposition 15
\begin{proposition}\label{11}
Let $(G_*,d^0_0,G)$ be an augmented pseudosimplicial group and
suppose that $d^0_0$ induces a natural isomorphism
$\pi_0(G_*)\ovs{\wt{d^0_0}}{\to}G$.
\begin{enumerate}
\item[(i)] Then the $n$-cube of groups ${\mathcal G}(n)$, $n\ge
1$, is induced by the normal $(n+1)$-ad of groups $(G_{n-1}, \Ker
d^{n-1}_0,\ldots, \Ker d^{n-1}_{n-1})$ i.e.
$$
{\mathcal G}(n)_A\cong G_{n-1}\diagup \uns{i\in A}{\prod} \Ker
d^{n-1}_{i-1}\;,\quad A\seq \la n\ra\;.
$$
\item[(ii)] $(G_*,
d^0_0, G)$ is aspherical if and only if the $n$-cube of groups
${\mathcal G}(n)$ is an exact $n$-presentation of the group $G$
for all $n\ge 1$.
\end{enumerate}
\end{proposition}
\begin{proof} {\bf (i)} Straightforward from the following
well-known fact on pseudosimplicial groups:
$$
d^n_j(\Ker d^n_i)=\Ker d^{n-1}_i\;\;\text{for}\;\; n>0\;,\; 0\le i
< j\le n\;.
$$
{\bf (ii)} It is well-known that asphericity of the augmented
pseudosimplicial group $(G_*,d^0_0, G)$ is equivalent to the
simplicial exactness of $(G_*,d^0_0,G)$, which certainly is
equivalent to the fact that ${\mathcal G}(n)$ is an exact
$n$-presentation of $G$ for all $n\ge 1$.
\end{proof}

\medskip

%%Remark
{\noindent \em Remark.} From Proposition \ref{11}(i) follows that
if $(F_*,d^0_0,G)$ is an augmented pseudosimplicial group such
that $d^0_0$ induces a natural isomorphism
$$\pi_0(F_*)\ovs{\wt{d^0_0}}{\to}G$$ and $F_i$, $i\ge 0$ are  free
groups, then the normal $(n+1)$-ad of groups $$(F_{n-1}, \Ker
d^{n-1}_0,\ldots ,\Ker d^{n-1}_{n-1})$$ satisfies the conditions
of Theorem BE (see also \cite{be}, \cite{e}). Thus, if these
conditions \emph{were} sufficient, simplicial exactness (or
asphericity) of $(F_*,d^0_0,G)$ would  not be needed for getting
the generalised Hopf formulas for higher homology of groups. It
was this, in fact, that made us suspect that these BE conditions
are not sufficient.

\medskip

Given an $n$-cube of groups ${\mathfrak F}$, a normal $(n+1)$-ad
of groups $(F;R_1,\ldots,R_n)$,  where $F={\mathfrak F}_{\eset}$
and $R_i=\Ker \al^{\eset}_{\{i\}}$, $i\in \la n\ra$ is called the
\emph{normal $(n+1)$-ad  of groups induced by ${\mathfrak F}$}. If
${\mathfrak F}$ is an exact $n$-presentation of the group
${\mathfrak F}_{\la n\ra}$, then the normal $(n+1)$-ad of groups
$(F;R_1,\ldots, R_n)$ satisfies the following condition:
$$
F\diagup \uns{i\in A}{\prod}R_i \cong {\mathfrak F}_A
\;\;\text{for all}\;\; A\seq \la n\ra\;,
$$
i.e. the $n$-cube of groups ${\mathfrak F}$ is induced by the normal
$(n+1)$-ad of groups $(F;R_1,\ldots, R_n)$ (see \cite{K2}).

\medskip

Now let ${\mathfrak F}$ be an $n$-presentation of the group $G$.
Applying $\check{C}$ (see Section 1) in the $n$-independent
directions, this construction leads naturally to an augmented
$n$-simplicial group. Taking the diagonal of this augmented
$n$-simplicial group gives the augmented simplicial group
$(\check{C}^{(n)}({\mathfrak F})_*,\al,G)$ called an
\emph{augmented $n$-fold \v{C}ech complex for} ${\mathfrak F}$,
where $\al=\al^{\eset}_{\la n\ra}:{\mathfrak F}_{\eset}\to G$. In
case ${\mathfrak F}$ is a free exact $n$-presentation of the group
$G$, then $(\check{C}^{(n)}({\mathfrak F})_*,\al,G)$ will be
called an \emph{$n$-fold \v{C}ech resolution of }$G$.

Let $G$, $H$ be groups. Let ${\mathfrak F}$ and $\mathfrak Q$ be
$n$-presentations of $G$ and $H$ respectively and $\lam:G\to H$ a
morphism of groups. A morphism $\ka:{\mathfrak F}\to {\mathfrak
Q}$ of $n$-cubes will be called an \emph{extension} of the group
morphism $\lam$ if $\ka_{\la n\ra}=\lam$.

\medskip

%%Theorem 16
\begin{theorem}\label{12}
Let ${\mathfrak F}$ and ${\mathfrak Q}$ be free and exact
$n$-presentations of given groups $G$ and $H$ respectively. Then
any morphism of groups $\lam:G\to H$ can be extended to a morphism
$\ka :{\mathfrak F}\to {\mathfrak Q}$ of $n$-cubes of groups which
naturally induces a morphism $\wt\ka$ of simplicial groups
$$
\begin{matrix}
\check{C}^{(n)}({\mathfrak F})_* & \ovs{\al}{\to} & G\\
\vcenter{\llap{$_{\wt\ka}{}$}}\da & & \da\vcenter{\rlap{${}_{\lam}$}} \\
\check{C}^{(n)}({\mathfrak Q})_* & \ovs{\al}{\to} & H
\end{matrix}
$$
over $\lam$. Furthermore, any two such extensions  $\ka,
\pi:{\mathfrak F}\to {\mathfrak Q}$ of $\lam$ induce simplicially
homotopic morphisms $\wt{\ka}$, $\wt{\pi}$ of simplicial groups,
denoted by $\wt{\ka}\simeq \wt{\pi}$.
\end{theorem}
\begin{proof} We begin by showing the existence of a
morphism of $n$-cubes of groups $\ka:{\mathfrak F}\to {\mathfrak
Q}$ extending the morphism of groups $\lam:G\to H$.

Since ${\mathfrak F}$ is free and ${\mathfrak Q}$ is exact, there
exists a homomorphism $\ka_{\la n\ra \stmin \{i\}}:{\mathfrak
F}_{\la n\ra \stmin \{i\}}\to {\mathfrak Q}_{\la n\ra \stmin
\{i\}}$ for all $i\in \la n\ra$, such that $\al^{\la n\ra\stmin
\{i\}}_{\la n\ra}\ka_{\la n\ra\stmin \{i\}}=\lam \al^{\la
n\ra\stmin \{i\}}_{\la n\ra}$. Suppose for some $A\seq \la n\ra$
and for all $B\supset  A$, $B\seq \la n\ra$, there exists a
homomorphism $\ka_B:{\mathfrak F}_B\to {\mathfrak Q}_B$ such that
$\al^B_C \ka_B=\ka_C \al^B_C$, $C\esq B$. Then as an immediate
consequence one has the induced homomorphism
$\ol{\ka}:\uns{B\supset A}{\lim}{\mathfrak F}_B\to \uns{B\supset
A}{\lim}{\mathfrak Q}_B$. Using again the facts that ${\mathfrak
F}$ is free and ${\mathfrak Q}$ is exact there exists a
homomorphism $\ka_A:{\mathfrak F}_A\to {\mathfrak Q}_A$ such that
$\al_A\ka_A=\ol{\ka}\al_A$. Clearly the constructed morphism of
$n$-cubes of groups $\ka:{\mathfrak F}\to {\mathfrak Q}$ naturally
induces a unique morphism of augmented $n$-simplicial groups and
applying the diagonal gives a morphism of simplicial groups
$\wt{\ka}:\check{C}^{(n)}({\mathfrak F})_*\to
\check{C}^{(n)}({\mathfrak Q})_*$ over the morphism $\lam$.

We need to prove the remaining part of the assertion first in a
particular case.

\noindent{\bf Particular Case}. \emph{Let $\ka,\;\pi:{\mathfrak
F}\to {\mathfrak Q}$ be two extensions of the group  morphism
$\lam:G\to H$ and $l\in \la n\ra$. Let
$\ka^{\{l\}}=\pi^{\{l\}}:{\mathfrak F}^{\{l\}}\to {\mathfrak
Q}^{\{l\}}$, then the respective induced morphisms of simplicial
groups $\wt{\ka}, \wt{\pi}:\check{C}^{(n)}({\mathfrak F})_*\to
\check{C}^{(n)}({\mathfrak Q})_*$ over $\lam$ are simplicially
homotopic.}

The construction of $\check{C}^{(n)}$ directly implies that for
any $n$-cube of groups ${\mathfrak F}$,
$\check{C}^{(n)}({\mathfrak F})_*$  is the diagonal of a
bisimplicial group $F_{**}$ induced by applying  the ordinary
\v{C}ech complex construction $\check{C}$ \;to the morphism  of
simplicial groups $\check{C}^{(n-1)}({\mathfrak
F}^{\ol{\{l\}}})_*\to \check{C}^{(n-1)}({\mathfrak F}^{\{l\}})_*$.

By assumption the extensions $\ka$ and $\pi$ of the group morphism
$\lam$ induce a commutative diagram of simplicial groups
$$
\begin{matrix}
\check{C}^{(n-1)}({\mathfrak F}^{\ol{\{l\}}})_* &
\atwo{\wt{\ka'}}{\wt{\pi'}} &
\check{C}^{(n-1)}({\mathfrak Q}^{\ol{\{l\}}})_*\\
\da & & \da \\
\check{C}^{(n-1)}({\mathfrak F}^{\{l\}})_* &
\atwo{\wt{\ka''}}{\wt{\pi''}} & \check{C}^{(n-1)}({\mathfrak
Q}^{\{l\}})_*
\end{matrix}
\;\;\;,
$$
where $\wt{\ka''}=\wt{\pi''}$, which implies there are morphisms of
simplicial objects of simplicial groups
$$
\begin{matrix}
F_{**} & \atwo{\ol\ka}{\ol\pi} & Q_{**}\\
\da & & \da\\
\check{C}^{(n-1)}({\mathfrak F}^{\{l\}})_* &
\ovs{\wt{\ka''}=\wt{\pi''}}{\to} & \check{C}^{(n-1)}({\mathfrak
Q}^{\{l\}})_*
\end{matrix}
$$
over the morphism of simplicial groups $\wt{\ka''}=\wt{\pi''}$.

The following lemmas will be needed.

\medskip

%%Lemma 17
\begin{lemma}\label{13}
Let $G_{**}$, $H_{**}$ be bisimplicial groups and
$\al,\;\be:G_{**}\to H_{**}$ morphisms of bisimplicial groups. Let
there exist a vertical (horizontal) simplicial homotopy
$h^v$($h^h$) between the induced morphisms of simplicial groups
$\al_m,\;\be_m:G_{m*}\to H_{m*}$($G_{*m}\to H_{*m}$) for all $m\ge
0$, such that the following conditions hold: $$d^h_j h^v_i= h^v_i
d^h_j\hspace{2cm}(d^v_j h^h_i= h^h_i d^v_j).$$ Then the induced
morphisms of simplicial groups $\wt\al,\;\wt\be:\Delta G_*\to
\Delta H_*$ are simplicially homotopic, $\wt\al \simeq \wt\be$,
where $\Delta G_*$ and $\Delta H_*$ are the diagonal simplicial
groups of $G_{**}$ and $H_{**}$ respectively.
\end{lemma}
\begin{proof} We can construct the required homotopy in
the following way: \\
$h'_i=h^v_i s^h_i: G_{nn}\to H_{n+1,n+1}$,
$0\le i\le n$.

Now we have to check the standard identities for simplicial homotopy. In fact,
$$
\begin{array}{rcl}
d^v_0 d^h_0 h^v_0 s^h_0\!&=&\!d^v_0 h^v_0 d^h_0 s^h_0= d^v_0 h^v_0= \al_{nn}\;,\\
d^v_{n+1} d^h_{n+1} h^v_n s^h_n\!&=&\!d^v_{n+1} h^v_n d^h_{n+1} s^h_n= d^v_{n+1} h^v_n= \be_{nn}\;,\\
d^v_i d^h_i h^v_j s^h_j\!&=&\!d^v_i h^v_j d^h_i s^h_j = \left\lbrace
  \begin{array}{cl}
h^v_{j-1} d^v_i s^h_{j-1} d^h_i=h^v_{j-1} s^h_{j-1} d^v_i d^h_i\;,  i<j\\
h^v_j d^v_{i-1} s^h_j d^h_{i-1}=h^v_j s^h_j d^v_{i-1} d^h_{i-1}\;,  i>j+1
 \end{array} \right.\;,
\\
d^v_{j+1} d^h_{j+1} h^v_{j+1} s^h_{j+1}\!&=&\!d^v_{j+1} h^v_{j+1}
d^h_{j+1} s^h_{j+1}= d^v_{j+1} h^v_{j+1}= d^v_{j+1} h^v_j
d^h_{j+1} s^h_j= d^v_{j+1} d^h_{j+1} h^v_j s^h_j\;.
\end{array}
$$
\end{proof}

\medskip

%%Lemma 18
\begin{lemma}\label{14}
Let $(\check{C}(\al)_*,\al,G)$, $(\check{C}(\be)_*,\be,H)$ be
augmented \v{C}ech complexes and $f,g: \check{C}(\al)_* \to
\check{C}(\be)_*$ morphisms of simplicial groups over a given
group morphism $\lam: G\to H$. Then $f$ and $g$ are simplicially
homotopic, $f \simeq g$.
\end{lemma}
\begin{proof} We only construct the simplicial homotopy
and leave the checking of the corresponding identities to the
reader. We define $h_i:\check{C}(\al)_n \to \check{C}(\be)_{n+1}$,
$0\le i\le n$, by $$h_i(x_0,\ldots, x_n)=(g(x_0),\ldots, g(x_i),
f(x_i),\ldots, f(x_n))$$ for all $(x_0,\ldots, x_n)\in
\check{C}(\al)_n$.
\end{proof}

\medskip

Returning to the main proof, using Lemma \ref{14}, it is easy to
see that there exists a vertical homotopy $h^v$ between the
induced morphisms of simplicial groups
$\ol{\ka_m},\;\ol{\pi_m}:F_{m*}\to Q_{m*}$ for all $m\ge 0$, such
that $d^h_j h^v_i=h^v_i d^h_j$. Applying Lemma \ref{13} there is a
simplicial homotopy between the morphisms of simplicial groups
$\wt{\ka},\;\wt{\pi}:\check{C}^{(n)}({\mathfrak F})_*\to
\check{C}^{(n)}({\mathfrak Q})_*$.

Now we return to the {\bf general case}, showing for any two
extensions $\ka,\;\pi:{\mathfrak F}\to {\mathfrak Q}$ of a group morphism
$\lam: G\to H$ the existence of extensions $\ka_1,\ldots,
\ka_{n-1}:{\mathfrak F}\to {\mathfrak Q}$ of $\lam$ such that
$\wt{\ka}\simeq \wt{\ka_1}$, $\wt{\ka_1}\simeq \wt{\ka_2}$,
\ldots, $\wt{\ka_{n-2}}\simeq \wt{\ka_{n-1}}$,
$\wt{\ka_{n-1}}\simeq \wt{\pi}$ which, of course, implies that
$\wt{\ka}\simeq \wt{\pi}$. In fact, we can construct an
extension $\ka_1:{\mathfrak F}\to {\mathfrak Q}$ in the following way: let
$\ka^{\{1\}}_1=\ka^{\{1\}}$ and $\ka^{\la n\ra \stmin
\{1\}}_1=\pi^{\la n\ra \stmin \{1\}}$. We complete the
construction of $\ka_1$ using the technique above and the facts
that ${\mathfrak F}$ is a free and ${\mathfrak Q}$ is an exact
$n$-presentation of the groups $G$ and $H$ respectively.

We construct an  extension  $\ka_i$ for all $2\le i\le n-1$ as
follows: let $\ka^{\{i\}}_i=\ka^{\{i\}}_{i-1}$ and $\ka^{\la
n\ra \stmin \la i\ra}_i=\pi^{\la n\ra \stmin \la i\ra}$. We
complete again the constructing of $\ka_i$ using the
above technique and  the facts that ${\mathfrak F}$ is free
and ${\mathfrak Q}$ is exact.

The construction of $\ka_i$, $1\le i\le n-1$, and our already
proved particular case imply that $\wt{\ka}\simeq \wt{\ka_1}$,
$\wt{\ka_1}\simeq \wt{\ka_2}$, \ldots, $\wt{\ka_{n-2}}\simeq
\wt{\ka_{n-1}}$, $\wt{\ka_{n-1}}\simeq \wt{\pi}$.
\end{proof}

\medskip

Using this comparison theorem  we make the following

\medskip

%%Definition
{\noindent \em Definition.} Let $T:{\mathfrak Gr}\to {\mathfrak Gr}$ be a
covariant functor. Define \emph{$i$-th $n$-fold \v{C}ech derived functor}
${\mathcal L}_i^{n-\text{fold}}T:{\mathfrak Gr}\to {\mathfrak Gr}$, $i\ge
0$, of the functor $T$ by choosing for each $G$ in ${\mathfrak Gr}$, a
free exact $n$-presentation ${\mathfrak F}$ and setting
$$
{\mathcal L}_i^{n-\text{fold}}T(G)=
\pi_i(T\check{C}^{(n)}({\mathfrak F})_*)\;,
$$
where $(\check{C}^{(n)}({\mathfrak F})_*,\al,G)$ is the $n$-fold
\v{C}ech resolution of the group $G$ for the free exact
$n$-presentation ${\mathfrak F}$ of $G$.

\medskip

The $n$-fold \v{C}ech complexes and hence the $n$-fold \v{C}ech
derived functors are closely related to the diagonal of the
$n$-simplicial multinerve of crossed $n$-cubes of groups .
In particular, we have the following

\medskip

%%Lemma 19
\begin{lemma}\label{15}
Let ${\mathfrak F}$ be an $n$-presentation of a group $G$. There
is an isomorphism of simplicial groups
$$
\check{C}^{(n)}({\mathfrak F})_* \cong E^{(n)}({\mathcal M})_*\;,
$$
where ${\mathcal M}$ is the inclusion crossed $n$-cube of groups
given by the normal $(n+1)$-ad of groups $(F;R_1,\ldots,R_n)$
induced by ${\mathfrak F}$.
\end{lemma}
\begin{proof} For $n=1$, is done in Lemma \ref{1.1} and hence for general
$n$, both constructions gives an isomorphism of $n$-simplicial
groups. Applying the diagonal clearly gives the result.
\end{proof}

\medskip

The following theorem gives the  $n$-th $n$-fold \v{C}ech
derived functor of the functor $Z_k:{\mathfrak Gr}\to {\mathfrak Gr}$,
$k\ge 2$.

\medskip

%%Theorem 20
\begin{theorem}\label{16}
Let $G$ be a group and $k\ge 2$. Then there is an isomorphism
$$
{\mathcal L}_n^{n-\text{fold}}Z_k(G)\cong \frac{\uns{i\in \la n\ra
}{\cap}R_i \cap \Gm_k(F)} {D_k (F;R_1,\ldots,R_n)}\;,\;\;\;n\ge
1\;,
$$
where $(F;R_1, \ldots, R_n)$ is the normal $(n+1)$-ad of groups
induced by some free exact $n$-presentation ${\mathfrak F}$ of the
group $G$.
\end{theorem}
\begin{proof} By Definition and Lemma \ref{15}, ${\mathcal
L}_n^{n-\text{fold}}Z_k(G)\cong \pi_n( Z_k E^{(n)}({\mathcal
M}))_*$, where ${\mathcal M}$ is the inclusion crossed $n$-cube of
groups given by the normal $(n+1)$-ad of groups $(F; R_1, \ldots,
R_n)$. Hence using Proposition \ref{8} one has an isomorphism
${\mathcal L}_n^{n-\text{fold}}Z_k(G)\cong \pi_n(E^{(n)}{\mathcal
B}_k({\mathcal M})_*)$. Then, by Proposition \ref{10},
\begin{align}\label{B}
{\mathcal L}_n^{n-\text{fold}}Z_k(G)\cong \uns{l\in \la
n\ra}{\cap} \Ker({\mathcal B}_k({\mathcal M})_{\la n
\ra}\ovs{\wt{\mu_{l,\la n \ra}}}{\lra} {\mathcal B}_k({\mathcal
M})_{\la n \ra \stmin \{l\}})\;.
\end{align}

Now we set up the  inductive hypothesis. Let $n=1$, then
$${\mathcal L}_1^{1-\text{fold}}Z_k(G)\cong
\Ker\Big(\frac{R_1}{D_k(F;R_1)}\lra
\frac{F}{\Gm_k(F)}\Big)=\frac{R_1\cap \Gm_k(F)}{D_k(F;R_1)}.$$

Proceeding by induction, we suppose that the result is true for
$n-1$ and we will prove it for $n$.

Let us consider $l\in\la n\ra$. It is easy to check that ${\mathfrak
F}^{\ol{\{l\}}}$ is a free exact $(n-1)$-presentation of the free
group ${\mathfrak F}_{\la n\ra\stmin \{l\}}$. Here we have to use the
fact that if $G$ is a free group, then ${\mathcal
L}_i^{n-\text{fold}}T(G)=0$, $i>0$ and ${\mathcal
L}_0^{n-\text{fold}}T(G)\cong T(G)$ for any functor $T:{\mathfrak
Gr}\to {\mathfrak Gr}$. Thus,  because of our inductive hypothesis,
\begin{equation}\label{C}
{\mathcal L}_{n-1}^{(n-1)-\text{fold}}Z_k({\mathfrak F}_{\la
n\ra\stmin \{l\}})\cong \frac{\uns{i\in \la n\ra \stmin
\{l\}}{\cap}R_i \cap \Gm_k(F)} {D_k
(F;R_1,\ldots,R_{l-1},R_{l+1},\ldots,R_n)}=0 \;.
\end{equation}

Now from (\ref{B}) and (\ref{C}) one can easily deduce that there
is the isomorphism
$$
{\mathcal L}_n^{n-\text{fold}}Z_k(G)\cong \frac{\uns{i\in \la
n\ra}{\cap}R_i \cap \Gm_k(F)} {D_k(F;R_1,\ldots,R_n)}\;.
$$
\end{proof}

\medskip

Now we are ready to improve on Theorem BE, and moreover to express by
generalised Hopf formulae not only the non-abelian derived
functors of the functor $Z_2$, but also the derived functors  of
the  functors $Z_k$, $k\ge 2$.

\medskip

%%Theorem 21
\begin{theorem}\label{17}
Let $G$ be a group,  ${\mathfrak F}$ be a free exact
$n$-presentation of $G$ and $k\ge 2$.  Then
$$
L_n Z_k(G)\cong \frac{\uns{i\in \la n\ra }{\cap}R_i \cap
\Gm_k(F)}{D_k (F;R_1,\ldots,R_n)}\;,\;\;\;n\ge 1\;,
$$
where $(F;R_1,\ldots,R_n)$ is the normal $(n+1)$-ad of groups
induced by ${\mathfrak F}$.
\end{theorem}
\begin{proof} This directly follows from Corollary \ref{6},
Proposition \ref{11}(ii) and Theorem \ref{16}.
\end{proof}

\medskip

%%Remark
{\noindent \em Remark.} One can prove an analogous result in a
more general context of Theorem BE and Corollary 4. In particular,
for a given group $G$ and an exact $n$-presentation ${\mathfrak
F}$ of $G$ such that
$$
L_1Z_k({\mathfrak F}_{\eset})=0\;, \;\;L_rZ_k({\mathfrak F}_A)=0
\;\; \text{for}\;\;r=|A|\;,\;r=|A|+1
$$
with $A$ a non-empty proper subset of $\la n\ra$, there is an
isomorphism
$$
L_nZ_k(G)\cong \frac{\uns{i\in \la n\ra }{\cap}R_i \cap
\Gm_k(F)}{D_k (F;R_1,\ldots,R_n)}\;\; \text{for}\;\;n\ge 1\;,
$$
where $(F;R_1,\ldots,R_n)$ is the normal $(n+1)$-ad of groups
induced by ${\mathfrak F}$.

\medskip

Now we  concentrate our attention on the computation of $2$-fold
\v{C}ech derived functors of the functor $Z_k:{\mathfrak Gr}\to
{\mathfrak Gr}$, $k\ge 2$. Using the fact that $Z_k$, $k\ge 2$ is
a right exact functor one easily shows that ${\mathcal
L}_0^{2-\text{fold}}Z_k\cong Z_k$. Moreover, by Proposition
\ref{8}, Proposition \ref{10} and Lemma \ref{15}, ${\mathcal
L}_i^{2-\text{fold}}Z_k = 0$ for $i\ge 3$. To take into account
Theorem \ref{16}, it  only remains to compute the first $2$-fold
\v{C}ech derived functor of the functor $Z_k$.

\medskip

%%Lemma 22
\begin{lemma}\label{18}(cf. Conduch\'{e} \cite{dc2} and also, \cite{atmp5})
Let $${\mathcal M} = \left\{\;\;
\begin{matrix}
L & \ovs{\lam}{\to} & M \\
\vcenter{\llap{$_{\lam'}{}$}}\da & & \da\vcenter{\rlap{${}_{\mu}$}} \\
N & \ovs{\nu}{\to} & P
\end{matrix}\;\;
\right\}$$
 be a crossed square. Then
\begin{eqnarray*}H_0(C_*({\mathcal
M}))&=&P\diagup \Im \mu \Im\nu,\\ H_1(C_*({\mathcal M}))&\cong&
M\times _P N\diagup \Im\ka,\\ H_2(C_*({\mathcal M}))&=&\Ker\lam \cap
\Ker\lam',
\end{eqnarray*}
 where $C_*({\mathcal M})$ is the mapping cone complex
of groups $$L \ovs{\al}{\to} M\rtimes N \ovs{\be}{\to} P$$ with
$\al(l)=({\lam(l)}^{-1},\lam'(l))$, $\be(m,n)=\mu(m)\nu(n)$ for
all $l\in L$, $(m,n)\in M\rtimes N$ and $\ka$ is the natural
homomorphism from $L$ to $M\times_P N$.
\end{lemma}
\begin{proof} We only prove that $H_1(C_*({\mathcal
M}))\cong M\times _P N\diagup \Im\ka$. It is easy to check that
$f:\Ker \be\to M\times _PN$, given by $f(m,n)=(m^{-1},n)$ for all
$(m,n)\in \Ker \be$, is an isomorphism and $\Im f\al= \Im \ka$.
The other results are as easy as this part  to check.
\end{proof}

\medskip

%%Proposition 23
\begin{proposition}
For a given group $G$ and $k\ge 2$ there are isomorphisms of
groups
$$
L_1Z_k(G)\cong {\mathcal L}_1^{2-\text{fold}}Z_k(G)\cong
\frac{R_1 \Gm_k(F) \cap R_2 \Gm_k(F)}{(R_1\cap R_2)\Gm_k(F)}\;,
$$
where $(F;R_1, R_2)$ is a normal $3$-ad of groups induced by some free
exact $2$-presentation ${\frak F}$ of $G$.
\end{proposition}
\begin{proof} Begin with the first isomorphism. By the
construction, $\check{C}^{(2)}({\frak F})_*$ is the diagonal of a
bisimplicial group $F_{**}$ induced by applying the ordinary
\v{C}ech complex construction $\check{C}$ to the morphism of
\v{C}ech complexes $\check{C}({\frak F}^{\ol{\{1\}}})_*\to
\check{C}({\frak F}^{\{1\}})_*$. Now applying the right exact
functor $Z_k$ dimension-wise, denote the resulted bisimplicial
group $Z_k(F_{**})$. By \cite{Qu} there is a spectral sequence
$$
E^2_{pq}\;\Longrightarrow\;{\mathcal L}_{p+q}^{2-\text{fold}}Z_k(G)\;,
$$
where $E^2_{0q}=0$, $q>0$ and $E^2_{p0}\cong {\mathcal
L}_p^{1-\text{fold}}Z_k(G)$, $p\ge 0$. Hence there is an isomorphism
${\mathcal L}_1^{2-\text{fold}}Z_k(G)\ovs{\cong}{\to} {\mathcal
L}_1^{1-\text{fold}}Z_k(G)$. Now the required isomorphism follows from \cite{Pi},
(see also \cite{HI}).

Again use of Proposition \ref{8}, Proposition \ref{10} and Lemma
\ref{15} implies that
$$
{\mathcal L}_1^{2-\text{fold}}Z_k(G)\cong H_1(C_*({\mathcal
B}_k({\mathcal M}))\;,
$$
where ${\mathcal M}$ is the inclusion crossed square induced by
the normal $3$-ad of groups $(F;R_1,R_2)$. Since the homomorphisms
$$
\wt{\mu_1}:{\mathcal B}_k({\mathcal M})_{\{1\}}\to {\mathcal
B}_k({\mathcal M})_{\eset}
$$
and
$$
\wt{\mu_2}:{\mathcal
B}_k({\mathcal M})_{\{2\}}\to {\mathcal B}_k({\mathcal
M})_{\eset}
$$
are injections, using Lemma \ref{18} implies the second
isomorphism.
\end{proof}

\medskip

Note that for group-abelianization functor ${\it Ab=Z_2}$ we have
the following new interpretation of the second integral group homology
$$
H_2(G)\cong {\mathcal L}_1^{2-\text{fold}}{\it Ab}(G)\cong
\frac{R_1 [F,F] \cap R_2 [F,F]}{(R_1\cap R_2)[F,F]}\;.
$$

It is an interesting problem to investigate and to compute the
functors ${\mathcal L}_i^{n-\text{fold}}Z_k$ for $0< i< n$, $n\ge
3$.

\

%%Section 6
\section{Hopf type formulas in algebraic $K$-theory}

\

Let us recall the well-known definition of
$\uns{\uns{j}{\lla}}{\lim}^{(1)}$, the first derived functor of
the functor $\uns{\uns{j}{\lla}}{\lim}$ (inverse limit in the
category of groups)(see e.g. \cite{HI}). Let $\{A_j,\; p^j_{j+1}\}_j$ be
a countable inverse system of groups, then
$$
\uns{\uns{j}{\lla}}{\lim}^{(1)}\{A_j,\; p^j_{j+1}\}=\uns{j}{\prod}A_j\diagup \sim \;,
$$
where $\sim$ is an equivalence relation on the set
$\uns{j}{\prod}A_j$ defined as follows: $\{a_j\}\sim \{a'_j\}$ if
there exists $\{h_j\}$ such that
$\{h_j\}\{a_j\}\{p^j_{j+1}(h^{-1}_{j+1})\}=\{a'_j\}$.

Let $\{G^j_*,\; \psi^j_{j+1}\}_j$ be a countable inverse system of
pseudosimplicial groups $G^j_*$ with $\psi^j_{j+1}:G^{j+1}_*\to
G^j_*$ a fibration for all $j\ge 1$. Let
$G_*=\uns{\uns{j}{\lla}}{\lim}\{G^j_*,\; \psi^j_{j+1}\}$.

\medskip

%%Theorem 24
\begin{theorem}\label{20}\cite{HI}
There is a short exact sequence of groups
$$
0\lra \uns{\uns{j}{\lla}}{\lim}^{(1)}\pi _{n+1}(G^j_*)\lra
\pi_n(G_*)\lra \uns{\uns{j}{\lla}}{\lim}\pi_n(G^j_*)\lra 0
$$
for all $n\ge 0$.
\end{theorem}

\medskip

Let us define the functor $Z_{\infty}:{\mathfrak Gr}\to {\mathfrak Gr}$ as
follows: for a given group $G$,
$Z_{\infty}(G)=\uns{\uns{j}{\lla}}{\lim}Z_j(G)$; for a given group
homomorphism $\lam:G\to H$, $Z_{\infty}(\lam)$ is the group
homomorphism induced by the $Z_j(\lam)$.

It is known from \cite{K1} (see also \cite{HI}) that the
non-abelian left derived functors $L^{\mathcal P}_*Z_{\infty}$ of
the functor $Z_{\infty}:{\mathfrak Gr}\to {\mathfrak Gr}$ are
isomorphic to Quillen's $K$-functors. Thus using Theorem \ref{20}
we deduce that there is a short exact sequence of abelian groups
$$
0\lra \uns{\uns{j}{\lla}}{\lim}^{(1)}\pi _{n+1}(F^j_*(GL({\mathfrak
R})))\lra K_{n+1}({\mathfrak R})\lra
\uns{\uns{j}{\lla}}{\lim}~\pi_n(F^j_*(GL({\mathfrak R})))\lra
0\;,\;\;n\ge 0\;,
$$
where $F_*(GL({\mathfrak R}))\to GL({\mathfrak R})$ is a free
pseudosimplicial resolution of the group $GL({\mathfrak R})$ and
$F^j_*(GL({\mathfrak R}))=Z_j(F_*(GL({\mathfrak R})))$.

Now according to Corollary \ref{6} we obtain the following

\medskip

%%Theorem 25
\begin{theorem}\label{21}
Let ${\mathfrak R}$ be a ring with unit and
$(F_*,d^0_0,GL({\mathfrak R}))$ be a free pseudosimplicial
resolution of the general linear group $GL({\mathfrak R})$. Then
there is an exact sequence of abelian groups
\begin{align*}
& 0\lra \uns{\uns{j}{\lla}}{\lim}^{(1)}\left(\frac{(\uns{i\in \la
n+1\ra}{\cap}\Ker d^n_{i-1})\cap \Gm_j(F_n)}{D_j(F_n;\Ker
d^n_0,\ldots,\Ker d^n_n)}\right) \lra K_{n+1}\left({\mathfrak
R}\right)
\lra \\
& \lra \uns{\uns{j}{\lla}}{\lim}\left(\frac{(\uns{i\in \la
n\ra}{\cap}\Ker d^{n-1}_{i-1})\cap
\Gm_j(F_{n-1})}{D_j(F_{n-1};\Ker d^{n-1}_0,\ldots,\Ker
d^{n-1}_{n-1})}\right) \lra 0
\end{align*}
for $n\ge 1$.
\end{theorem}

\medskip

Note that using Theorem \ref{17} and its Remark, one can express
$K_{n+1}({\mathfrak R})$ in data coming from exact $(n+1)$ and
$n$-presentations of the group $GL({\mathfrak R})$. We hope to
return to a more detailed analysis of this in future work.

\

\

\begin{center}
{\bf Acknowledgements}
\end{center}
The second author would like to thank the Royal Society, and
University of Wales Bangor for their hospitality at the
`initialisation'  of this paper. The authors were partially
supported by INTAS grant  00 566, and in addition the third author
by INTAS 971 - 31961. The first and the second authors were also
supported by INTAS Georgia grant No 213, whilst the second author
was supported by FNRS grant 7GEPJ065513.01.

\

\end{document}